\magnification=\magstep1
\baselineskip = 14 pt
\font\bigbf=cmbx10 scaled \magstep1
\def\qn{{q^n}}
\def\qk{{q^k}}
\def\qnm{{q^{n-1}}}
\def\qkm{{q^{k-1}}}
\def\qmn{{q^{-n}}}
\def\qmnp{{q^{-n+1}}}
\def\qmkp{{q^{-k+1}}}
\def\qtn{{q^{2n}}}

\def\qtnp{{q^{2n+1}}}
\def\qntm{{q^{n-2}}}
\def\qnmt{{q^{n-2}}}
\def\qtnmt{{q^{2n-2}}}
\def\qmtn{{q^{-2n}}}
\def\qmtnp{{q^{-2n+1}}}

\def\qtnmr{{q^{2n-3}}}
\def\qmnpt{{q^{-n+2}}}
\def\qmkpt{{q^{-k+2}}}
\def\qnpt{{q^{n+2}}}

\def\qnp{{q^{n+1}}}
\def\qtnm{{q^{2n-1}}}
\def\qpmo{{q^{\pm 1}}}
\def\qtnpt{{q^{2n+2}}}
\def\qtmnpr{{q^{-2n+3}}}
\def\qpmt{{q^{\pm 2}}}
\def\qN{{q^N}}
\def\qmN{{q^{-N}}}
\def\qNpo{{q^{N+1}}}
\def\qnpo{{q^{n+1}}}
\def\qNpt{{q^{N+2}}}

\def\qtmN{{q^{2-N}}}
\def\qmk{{q^{-k}}}
\def\qkpo{{q^{k+1}}}
\def\qtkpo{{q^{2k+1} }}
\def\qtmn{{q^{2-n}}}
\def\qtnpo{{q^{2n+1} }}
\def\qtk{{q^{2k}}}
\def\qNpk{{q^{N+k}}}
\def\qmNpk{{q^{-N+k} }}
\def\qmNpn{{q^{-N+n} }}
\def\qomN{{q^{1-N} }}
\def\qmNpo{{q^{-N+1} }}
\def\qNmo{{q^{N-1} }}
\def\qNmt{{q^{N-2} }}
\def\qNmn{{q^{N-n} }}
\def\qmNpt{{q^{-N+2} }}
\def\qomn{{q^{1-n} }}
\def\qnmo{{q^{n-1}}}
\def\qtnmo{{q^{2n-1}}}
\def\qopn{{q^{1+n} }}
\def\qM{{q^M}}
\def\qmM{{q^{-M}}}
\def\qmm{{q^{-m}}}

\def\W{{_{10}W_9}}
\def\w{{_8W_7}}
\def\tphia{{_{10} \phi_9}}
\def\ephis{{_{8} \phi_7}}
\def\rphit{{_{3} \phi_2}}
\def\tphio{{_{2} \phi_1}}
\def\frac#1#2{{ {{#1} \over {#2}} }}
\def\coloneq{{\, :=\, }}
\def\Phib{{\Phi^{(b)} }}
\def\CW{{\cal W}}
\def\caprom#1{\rm{\uppercase\expandafter{\romannumeral #1}}}
\def\npo{{n+1}}
\def\nmo{{n-1}}
\def\lowminus {\ {{ {\vphantom{\lambda_2}} } \atop {
{\vphantom{\lambda_2}} -}}\ }
\def\lowdots {{{{\vphantom{\lambda_2}}} \atop {{\vphantom{
\lambda_2}} \dots}}}

\def\bigslash{{\slash}}

\centerline{\bigbf Contiguous Relations, Continued Fractions
and Orthogonality\footnote*{\rm Research partially supported
by NSERC (Canada).}}
\vglue 1 truein

\centerline{Dharma P.~Gupta}

\centerline{and}

\centerline{David R.~Masson}
\bigskip

\centerline{Department of Mathematics}

\centerline{University of Toronto, Toronto, M5A 1A1}
\vglue 1 truein

\centerline{\bf Abstract}
\medskip

We examine a special linear combination of balanced
very-well-poised $\tphia$ basic hypergeometric series
that is known to satisfy a transformation. We call this
$\Phi$ and show that it satisfies certain three-term
contiguous relations. From two sets of contiguous
relations for $\Phi$ we obtain fifty-six pairwise
linearly independent solutions to a three-term recurrence
that generalizes the recurrence for Askey-Wilson
polynomials. The associated continued fraction is
evaluated using Pincherle's theorem. From this continued
fraction we are able to derive a discrete system of
biorthogonal rational functions. This ties together
Wilson's results for rational biorthogonality, Watson's
$q$-analogue of Ramanujan's Entry 40 continued fraction
and a conjecture of Askey concerning the latter. Some new
$q$-series identities are also obtained. One is an
important three-term transformation for $\Phi$'s which
generalizes all the known two and three-term $\ephis$
transformations. Others are new and unexpected quadratic
identities for these very-well-poised $\ephis$'s.

\vfill



\eject

\noindent{\bf 1. Introduction}
\medskip
\noindent{\bf Background:} It is now ten years since Askey-Wilson 
polynomials  were introduced as an explicit system of orthogonal
polynomials that generalize Jacobi polynomials [2]. Here we are
concerned with a generalization of Askey-Wilson
polynomials which yields an explicit system of biorthogonal
rational functions. This generalization was first given by Wilson
in 1991 [30].

The explicitness of the above orthogonalities is connected to the
fact that certain hypergeometric and basic hypergeometric
series satisfy three-term contiguous relations. The
importance of contiguous relations in this context was
emphasized by Wilson in his 1978 thesis [29]. For more recent work
which makes extensive use of contiguous relations see [16], 
[23], [13], [14], [15].

In the case of Askey-Wilson polynomials, the relevant series
are a terminating, balanced ${}_4 \phi_3$ series, which
describes the polynomials themselves [2], and a very-well-poised
$\ephis$ series, which describes the minimal solution to the
polynomial recurrence [18], [11]. 

The necessity of a three-term contiguous relation is
dictated by the fact that an orthogonal polynomial system
$\{ P_n (x) \}_{n=0}^\infty$ must satisfy a three-term
recurrence of the form [6]
$$
P_{n+1} (x) - (x-c_{n+1} ) P_n (x) + \lambda_{n+1} P_{n-1}
(x) =0,
\leqno{(1.1)}
$$
associated with the $J$-fraction
$$
\frac{1}{x-c_1} \lowminus
\frac{\lambda_2}{x-c_2} \lowminus
\frac{\lambda_3}{x-c_3} \lowminus \lowdots . 
\leqno{(1.2)}
$$

Are there explicit orthogonal polynomial systems more
general than Askey-Wilson polynomials? The answer is
believed to be no since there are no known basic
hypergeometric series which are more general and which also
satisfy a three-term contiguous relation that can be cast
into the form (1.1).

However, if we relax the requirement that the orthogonal
system consists of polynomials, then we are led to Wilson's
system of biorthogonal rational functions [30]. With only
a slight modification this again seems to be the most general
model of its type.

In this paper we examine this general rational function 
biorthogonality by starting with some three-term contiguous 
relations of the most general known type. The basic 
hypergeometric series involved are terminating, balanced, 
very-well-poised ${}_{10}\phi_9$'s or, more generally, 
special linear combinations of non-terminating such 
${}_{10}\phi_9$'s which we call 
$\Phi$'s. Both of these satisfy three-term contiguous 
relations which can be put into the form
$$
P_{n+1} (x) - (x-c_{n+1})P_n (x) + \lambda_{n+1} (x -
\alpha_{n+1})(x -
\beta_{n+1} )P_{n-1} (x) =0
\leqno{(1.3)}
$$
associated with the $R_{II}$-fraction
$$
\frac{1}{x-c_1} \lowminus
\frac{\lambda_2(x-\alpha_2)(x-\beta_2)}{x-c_2} \lowminus
\frac{\lambda_3(x-\alpha_3)(x-\beta_3)}{x-c_3} \lowminus \lowdots
.
\leqno{(1.4)}
$$
That the orthogonal systems corresponding to (1.3) and
(1.4) involve rational functions rather than polynomials
has been demonstrated by Ismail and Masson [17].

There is a second related aspect which lead us to examine
this most general basic hypergeometric level. Many of
the continued fractions of Ramanujan are connected with
orthogonal polynomials. Some of the most intriguing of
these are expressed in terms of ratios of products of
gamma functions [5] and have been shown to be connected with
special and limiting cases of Askey-Wilson polynomials [21], [22].
Armed
with this fact we were able to extend and give new meaning to
Ramanujan's famous Entry 40 [24] and its $q$-analogue given by
Watson [28], [10]. These are all special cases of the very
general continued fraction we examine in Section 3 of
this paper. From the simplest terminating form of this continued 
fraction given in
Corollary 3.3 and the orthogonality derived in Section 4,
we are now able to vindicate Askey's conjecture [1, p. 37] that 
Ramanujan's Entry 40 is connected with Dougall's
${}_7F_6$ summation formula and a three-term recurrence for a 
very-well-poised ${}_9F_8$.
\medskip
\noindent{\bf Results:} In Section 2 we prove that $\Phi$
satisfies certain three-term 
contiguous relations. For two of these (Theorems 2.4
and 2.5) we make essential use of the two term
transformation formula for $\Phi$ [7, ($\caprom 3$.~39),
p.~247]. This generalizes our earlier work in [10].

In Section 3, Theorems 2.4 and 2.5 are used to obtain
fifty-six pairwise linearly independent solutions to a very
general eight parameter three-term recurrence. This
recurrence generalizes the recurrence for Askey-Wilson
polynomials. The large $n$ asymptotics of the solutions
is also examined. From two of these solutions we are
then able to construct a minimal solution. This minimal
solution, via Pincherle's Theorem, gives an explicit but
complicated continued fraction result which generalizes
Ramanujan's Entry 40 and our earlier work [10]. Two special
cases are considered. One is a terminating fraction given
in terms of a terminating ${}_{10}\phi_9$ [28].

In Section 4 this terminating fraction is expressed as
an $R_{II}$-fraction and used to derive a very general
explicit rational biorthogonality. This derivation
follows the methods in [17]. Also in Section 4 we give
six limiting cases of this biorthogonality. Five limits
are at the ${}_4\phi_3$ or $\ephis$ level and one is a
$q \rightarrow 1$ limit at the ${}_9F_8$ level. One of
the ${}_4\phi_3$ limits corresponds to the case of
$q$-Racah polynomials while the other orthogonalities
are new. By taking further limits this can be extended
to a full Askey-type scheme of rational biorthogonality.
For a detailed review of the Askey-scheme for polynomial
orthogonality see Koekoek and Swarttouw [20].

In Section 5 we give further details for one of the limit 
cases at the very-well-poised $\ephis$ level. This expands 
on some of our previous results for this model [12]. Also 
a general three-term $\ephis$ transformation is derived 
which is essential for Section 6. 

We have already mentioned that $\Phi$ satisfies a two term
transformation formula. This is the most general two term
transformation given in Gasper and Rahman [7].  In Section
6 we derive a missing companion transformation. This is
a three-term $\Phi$ transformation which generalizes all
the known $\ephis$ transformations.

Finally, in Section 7 we derive what we feel are some unexpected
$\ephis$
identities. These new quadratic identities are derived from the
asymptotics
of the Casorati determinants of some of the recurrence solutions of
Section 3.
\medskip

\noindent{\bf Notation:} We follow the notation in Gasper
and Rahman [7] except that we omit the designation `$q$'
for the base in the $q$-shifted factorials and the basic
hypergeometric functions. Thus, given a number $q$, the
$q$-shifted factorial is defined by
$$
\eqalign{
& (a)_0 \coloneq 1, \qquad (a)_n \coloneq
\prod\limits_{j=1}^n (1-aq^{j-1}), \quad n=1,2,\dots \cr
& (a)_{-n} \coloneq {1 \over {(aq^{-n})_n}}, \quad n=1,2,
\dots, \cr
}
$$
and for $\vert q \vert < 1$, 
$$
(a)_\infty = \prod\limits_{j=1}^\infty (1 - aq^{j-1}).
$$
We also write 
$$
(a_1, a_2, \dots, a_k)_n \coloneq \prod\limits_{j=1}^k
(a_j)_n, \quad n \ {\rm integer\ or}\ n=\infty.
$$
The $_r \phi_s$ basic hypergeometric series is
given by
$$
_r \phi_s \left( {{a_1,a_2,\dots, a_r} \atop
{b_1,b_2, \dots, b_s}}; z \right) \coloneq
\sum\limits_{n=0}^\infty {{(a_1,a_2,\dots, a_r)_n} \over
{(b_1, b_2, \dots, b_s, q)_n}} \left[ (-1)^n q^{\left( {n
\atop 2} \right)} \right]^{1 + s -r} z^n,
\leqno{(1.1)}
$$
and the general bilateral basic hypergeometric series is
defined by 
$$
_r \psi_s \left( {{a_1,a_2,\dots, a_r} \atop
{b_1,b_2, \dots, b_s}}; z \right) \coloneq
\sum\limits_{n=-\infty}^\infty {{(a_1,a_2,\dots, a_r)_n} \over
{(b_1, b_2, \dots, b_s)_n}} \left[ (-1)^n q^{\left( {n
\atop 2} \right)} \right]^{s -r} z^n.
\leqno{(1.2)}
$$
We denote a very-well-poised, balanced $\tphia$ series as 
$$
\eqalign{
& \phi = \W(a;b,c,d,e,f,g,h;q) \coloneq \tphia \left(
{{a, q\sqrt{a}, -q\sqrt{a},b,c,d,e,f,g,h} \atop {\sqrt{a},
-\sqrt{a}, \frac{aq}{b}, \frac{aq}{c}, \frac{aq}{d},
\frac{aq}{e},
\frac{aq}{f},
\frac{aq}{g},
\frac{aq}{h} }}; q \right), \cr
& \qquad bcdefgh=a^3q^2 ,\cr
}
\leqno{(1.3)}
$$
and the limiting case of a very-well-poised $\ephis$ as
$$ 
\w(a;b,c,d,e,f;{a^2q^2\over bcdef}) \coloneq _8\phi_7 \left(
{{a, q\sqrt{a}, -q\sqrt{a},b,c,d,e,f} \atop {\sqrt{a},
-\sqrt{a}, \frac{aq}{b}, \frac{aq}{c}, \frac{aq}{d}, \frac{aq}{e}, 
\frac{aq}{f}}}; {a^2q^2\over bcdef}\right). \leqno{(1.3')}
$$
A complementary pair of very-well-poised, balanced 
$\tphia$'s is defined by
$$
\eqalign{
& \Phi^{(b)} (a;b,c,d,e,f,g,h;q) \cr
& \qquad \coloneq
\W(a;b,c,d,e,f,g,h;q) + 
{{
\left(
aq,
\frac{b}{a},c,d,e,f,g,h, \frac{bq}{c}, 
\frac{bq}{d},
\frac{bq}{e},
\frac{bq}{f},
\frac{bq}{g},
\frac{bq}{h}\right)_\infty } \over {
\left( \frac{b^2q}{a}, \frac{a}{b}, 
\frac{aq}{c},
\frac{aq}{d},
\frac{aq}{e},
\frac{aq}{f},
\frac{aq}{g},
\frac{aq}{h},
\frac{bc}{a},
\frac{bd}{a},
\frac{be}{a},
\frac{bf}{a},
\frac{bg}{a},
\frac{bh}{a} \right)_\infty }} 
\cr
& \qquad\qquad  
\times \W \left( \frac{b^2}{a}; b, 
\frac{bc}{a},
\frac{bd}{a},
\frac{be}{a},
\frac{bf}{a},
\frac{bg}{a},
\frac{bh}{a}; q \right) = \phi + \phi'^{(b)}, {\rm\ say,\ with} 
\cr
& \qquad bcdefgh=a^3q^2. \cr}
\leqno{(1.4)}
$$
Here, `$b$' is a {\it distinguished parameter\/} which can
be interchanged with $c$, $d$, $e$, $f$, $g$ or $h$ to give
different $\Phi$'s. We use the notation 
$$
\Phi  = \Phi (a;b,c,d,e,f,g,h;q)=\phi+\phi'
\leqno{(1.5)}
$$
to denote any one of these seven possible complementary
pairs. Whenever we want to specify the distinguished
parameter `$b$' while defining $\Phi$ or $\phi'$, we shall
denote them by $\Phi^{(b)}$ and $\phi'^{(b)}$ respectively.

We will avoid singularities in the definition of
$\Phi^{(b)}$. In particular, we assume that $a \neq b$. In case any
of the
parameters $c$, $d$, $e$, $f$, $g$, $h$ is $1$, the
complementary part $\phi'^{(b)}$ vanishes.

A $q$-analogue of Wilson's notation [29] will be used for
the variations of $\phi$, $\Phi$ or $\phi'$ with respect
to the parameters.  Thus, $\Phib (g+, h-)$ represents the
expression (1.4) with `$g$' and `$h$' replaced throughout
by `$gq$' and `$\frac{h}{q}$' respectively. $\Phib_\pm$
denotes the expression which would be obtained by
replacements
$$
(a,b,c,d,e,f,g,h) \rightarrow (a\qpmt, b\qpmo,
c\qpmo, d\qpmo, e\qpmo, f\qpmo, g\qpmo, h\qpmo)
$$
throughout the expression (1.4).

\bigskip
\noindent{\bf 2. Contiguous Relations}
\medskip

We obtain three-term contiguous relations for the
$\Phi$ function which generalize the relations
we obtained earlier [10] for $\tphia$'s. The method of
proof is basically the same as the one employed by
Wilson [29] in his study of $_9F_8$ functions and
employed by us in our previous work [10].

\medskip
\noindent{\bf Lemma 2.1.} [10, p. 431, Lemma 1]
{\it If $\phi$ denotes the balanced very-well-poised
$\tphia$ series defined by (1.3), then}
$$
\eqalign{
& \phi (g-, h+) - \phi \cr
& \quad = 
{{ {{aq} \over h} 
(1 - {{hq} \over g}) (1 - {{gh} \over {aq}})
(1-aq) (1-aq^2) (1-b) (1-c) (1-d) (1-e) (1-f)} \over {
(1 - {{aq} \over {g}}) (1 - {{aq^2} \over {g}})
(1 - {{a} \over {h}}) (1 - {{aq} \over {h}})
(1 - {{aq} \over {b}}) (1 - {{aq} \over {c}})
(1 - {{aq} \over {d}}) (1 - {{aq} \over {e}})
(1 - {{aq} \over {f}}) }} \phi_{+} (g-). \cr
}
\leqno{(2.1)}
$$

\medskip
\noindent{\bf Lemma 2.2} 
{\it If $\phi$ denotes the balanced very-well-poised
$\tphia$ series defined by {\rm (1.3)}, then}
$$
\eqalign{
c\left(1-c\right)
\left(1 - \frac{a}{c}\right)
\left(1 - \frac{dq}{g}\right)
\left(1 - \frac{gd}{aq}\right)
\phi \left(g-, c+\right) & \cr
-d\left(1-d\right)
\left(1 - \frac{a}{d}\right)
\left(1 - \frac{cq}{g}\right)
\left(1 - \frac{gc}{aq}\right)
\phi \left(g-, d+\right) & \cr
+d\left(1-\frac{g}{q}\right)
\left(1 - \frac{c}{d}\right)
\left(1 - \frac{aq}{g}\right)
\left(1 - \frac{cd}{a}\right)
\phi & = 0. \cr
}
\leqno{(2.2)}
$$

\noindent{\bf Proof.} Elimination of $\phi_+ (g-)$ from
two relations written for $\phi (g-, c+) \break -\phi$ and
$\phi (g-,d+) -\phi$ with the help of Lemma 1 gives the
required result. 

\medskip
\noindent{\bf Lemma 2.1$'$.} {\it Irrespective of the choice
of the distinguished parameter, $\Phi$ 
satisfies the relation}
$$
\eqalign{
& \Phi (g-, h+) - \Phi \cr
& \quad = 
{{ {{aq} \over h} 
(1 - {{hq} \over g}) (1 - {{gh} \over {aq}})
(1-aq) (1-aq^2) (1-b) (1-c) (1-d) (1-e) (1-f)} \over {
(1 - {{aq} \over {g}}) (1 - {{aq^2} \over {g}})
(1 - {{a} \over {h}}) (1 - {{aq} \over {h}})
(1 - {{aq} \over {b}}) (1 - {{aq} \over {c}})
(1 - {{aq} \over {d}}) (1 - {{aq} \over {e}})
(1 - {{aq} \over {f}}) }} \Phi_{+} (g-). \cr
}
\leqno{(2.3)}
$$

\noindent{\bf Proof.} In view of Lemma 2.1, we only need
to prove that (2.3) holds true for the complementary
part $\phi'$ irrespective of the choice of the
distinguished parameter. It is evident that we should
check the validity of the statement in three different
cases viz., when the distinguished parameter is either
`$g$' or `$h$' or it is one of the parameters $b$, $c$,
$d$, $e$ or $f$.

The statement (2.3) can be shown to be true for 
$\phi'^{(g)}$ if we apply Lemma 1 to the series
$$
\W \left( \frac{g^2}{aq^2}; \frac{gb}{aq},
\frac{gc}{aq}, \frac{gd}{aq}, \frac{ge}{aq},
\frac{gf}{aq}, g, \frac{gh}{aq}; q \right)
\leqno{(2.4)}
$$
and make the required simplification. 

Similarly, validity of (2.3) for $\phi'^{(h)}$ can be
derived by applying Lemma 1 to
$$
\W \left( \frac{h^2}{a}; \frac{hb}{a}, \frac{hc}{a},
\frac{hd}{a}, \frac{he}{a}, \frac{hf}{a}, \frac{hg}{a},
h; q \right).
\leqno{(2.5)}
$$
In the third case, say for example for $\phi'^{(b)}$,
we can apply (2.2) of Lemma 2 to the function 
$$
\W \left( \frac{b^2}{a}; b, \frac{bc}{a},
\frac{bd}{a}, \frac{be}{a}, \frac{bf}{a}, \frac{bg}{a},
\frac{bh}{a}; q \right)
\leqno{(2.6)}
$$
and we arrive at the desired result. This completes the
proof of the Lemma. 

In view of Lemma 2.2 and Lemma 2.1$'$ we can immediately
state Lemma 2.2$'$.

\medskip
\noindent{\bf Lemma 2.2$'$.} 
{\it Irrespective of the choice of the distinguished parameter,
$\Phi$ satisfies the relation}
$$
\eqalign{
c\left(1-c\right)
\left(1 - \frac{a}{c}\right)
\left(1 - \frac{dq}{g}\right)
\left(1 - \frac{gd}{aq}\right)
\Phi \left(g-, c+\right) & \cr
-d\left(1-d\right)
\left(1 - \frac{a}{d}\right)
\left(1 - \frac{cq}{g}\right)
\left(1 - \frac{gc}{aq}\right)
\Phi \left(g-, d+\right) & \cr
+d\left(1-\frac{g}{q}\right)
\left(1 - \frac{c}{d}\right)
\left(1 - \frac{aq}{g}\right)
\left(1 - \frac{cd}{a}\right)
\Phi & = 0 \cr
}
\leqno{(2.7)}
$$

\medskip
\noindent{\bf Theorem 2.3.} {\it Irrespective of the choice of the
distinguished parameter, $\Phi$ satisfies the relation}
$$
\eqalign{
{{g^2 (1-h)
(1 - \frac{aq}{gb})
(1 - \frac{aq}{gc})
(1 - \frac{aq}{gd})
(1 - \frac{aq}{ge})
(1 - \frac{aq}{gf})
} \over {
(1- \frac{aq}{g})
(1- \frac{aq^2}{g}) 
}} \Phi_+ (g-) & \cr
- {{h^2 (1-g)
(1 - \frac{aq}{hb})
(1 - \frac{aq}{hc})
(1 - \frac{aq}{hd})
(1 - \frac{aq}{he})
(1 - \frac{aq}{hf})
} \over {
(1- \frac{aq}{h})
(1- \frac{aq^2}{h}) 
}} \Phi_+ (h-) & \cr
- {{g (1-{h \over g})
(1 - \frac{aq}{b})
(1 - \frac{aq}{c})
(1 - \frac{aq}{d})
(1 - \frac{aq}{e})
(1 - \frac{aq}{f})
} \over {
(1- {aq})
(1- {aq^2}) 
}} \Phi & =0. \cr
}
\leqno{(2.8)}
$$

\noindent{\bf Proof.}
First we  indicate the proof for 
$$
\Phib (a;b,c,d,e,f,g,h;q).
$$
Assuming that (2.8) is valid, we apply Bailey's transformation
(see Gasper and Rahman [7], (2.30), p. 56) to each of 
$\Phib_+ (g-)$, $\Phib_+ (h-)$ and $\Phib$.
Subsequently we have three pairs of $\tphia$'s
corresponding to the three $\Phi$'s. We can pick out
one $\tphia$ from each of the three pairs so that this
set of three $\tphia$'s and the remaining set of three
both separately satisfy valid
three-term contiguous relations. In fact, one set satisfies the
relation we obtain by applying Lemma 2.1 to 
$\W
\left( \frac{bh}{g}; b, \frac{aq}{cg}, \frac{aq}{dg},
\frac{aq}{eg}, \frac{aq}{fg}, \frac{bh}{a}, h; q \right)$.
The other set satisfies the relation obtained by
applying Lemma 2.1 to
$\W
\left( \frac{bg}{h}; b, \frac{aq}{ch}, \frac{aq}{dh},
\frac{aq}{eh}, \frac{aq}{fh}, \frac{bh}{a}, g; q \right)$.
This completes the proof of (2.8) when `$b$' is the
distinguished parameter. We have a similar proof when
`$h$' is the distinguished parameter. This time one
combination of three $\tphia$'s gets disposed of exactly
as in the previous case. The other combination of
three $\tphia$'s constitutes a relation which is the same
as that obtained by applying Lemma 2.2 to the series
$$
\W \left( \frac{hg}{b}; h, \frac{aq}{bc},
\frac{aq}{bd}, \frac{aq}{be}, \frac{aq}{bf},
\frac{hg}{aq}, gq; q \right).
$$
This completes the proof of Theorem 2.3.

\medskip
\noindent{\bf Theorem 2.4.} 
{\it Irrespective of the choice of the distinguished
parameter, $\Phi$ satisfies the relation}
$$
\eqalign{
{{g (1-h)
(1-\frac{a}{h})(1-\frac{aq}{h})
(1 - \frac{aq}{gb})
(1 - \frac{aq}{gc})
(1 - \frac{aq}{gd})
(1 - \frac{aq}{ge})
(1 - \frac{aq}{gf})
} \over {
(1- \frac{hq}{g})
}} \left[ \Phi (g-, h+) -\Phi \right] & \cr
- {{h (1-g)
(1-\frac{a}{g})(1-\frac{aq}{g})
(1 - \frac{aq}{hb})
(1 - \frac{aq}{hc})
(1 - \frac{aq}{hd})
(1 - \frac{aq}{he})
(1 - \frac{aq}{hf})
} \over {
(1- \frac{gq}{h})
}} \left[ \Phi (h-, g+) -\Phi \right] & \cr
- {{aq} \over h}
\left(1 - \frac{h}{g}\right)
\left(1 - \frac{gh}{aq}\right)
(1 - b)
(1 - c)
(1 - d)
(1 - e)
(1 - f)
\Phi & =0. \cr
}
\leqno{(2.9)}
$$

\noindent{\bf Proof.} 
The result can be obtained from Theorem 2.3, if we
substitute in (2.8) the value of $\Phi_+(g-)$ as given by
(2.3) and the value of $\Phi_+(h-)$ which would be
obtained from a $g \leftrightarrow h$ interchange of
(2.3).

\medskip
\noindent{\bf Theorem 2.5.}
{\it Irrespective of the choice of the distinguished
parameter, $\Phi$ satisfies the relation}
$$
\eqalign{
& {{agq} \over h} 
{{
(1-aq)(1-aq^2)
(1 -\frac{aq}{gb})
(1 -\frac{aq}{gc})
(1 -\frac{aq}{gd})
(1 -\frac{aq}{ge})
(1 -\frac{aq}{gf})
} \over {
(1 -\frac{aq}{g})
(1 -\frac{aq^2}{g})
(1 -\frac{aq}{b})
(1 -\frac{aq}{c})
(1 -\frac{aq}{d})
(1 -\frac{aq}{e})
(1 -\frac{aq}{f})
}} \cr
& \qquad \times \left(1-\frac{gh}{aq}\right) 
(1-h)
(1-b)
(1-c)
(1-d)
(1-e)
(1-f)
\Phi_+(g-) 
\cr
& \qquad - {{q
(1- \frac{a}{g})
(1- \frac{aq}{g})
(1- \frac{a}{h})
(1- \frac{aq}{h})
(1- \frac{a}{b})
(1- \frac{a}{c})
(1- \frac{a}{d})
(1- \frac{a}{e})
(1- \frac{a}{f})
} \over {
(1- \frac{a}{q})
(1-a)
}}
\Phi_- (g+) \cr
& \qquad - \left[ \vphantom{{{ \left({a \over b}\right)} \over
{\left({a \over b}\right)}}}
\frac{aq}{h} 
\left(1-\frac{h}{g}\right)
\left(1-\frac{gh}{aq}\right)
(1-b)
(1-c)
(1-d)
(1-e)
(1-f)
\right. \cr
& \qquad + \frac{g^2q^2}{h} {{
(1- \frac{aq}{g})
(1- \frac{a}{h})
(1- \frac{aq}{h})
(1- \frac{h}{q})
(1- \frac{a}{gb})
(1- \frac{a}{gc})
(1- \frac{a}{gd})
(1- \frac{a}{ge})
(1- \frac{a}{gf})
} \over {
(1- \frac{gq}{h})
(1- \frac{a}{gq})
}} \cr
& \left.  \qquad - {{
h (1-g)
(1-\frac{a}{g})
(1-\frac{aq}{g})
(1-\frac{aq}{hb})
(1-\frac{aq}{hc})
(1-\frac{aq}{hd})
(1-\frac{aq}{he})
(1-\frac{aq}{hf})
} \over { (1 - \frac{gq}{h})}} \right] \Phi = 0. \cr
}
\leqno{(2.10)}
$$

\noindent{\bf Proof.} 
First we make the parameter replacements $(a,b,c,d,e,f,g,h)
\rightarrow \break (\frac{a}{q^2}, \frac{b}{q},
\frac{c}{q}, \frac{d}{q}, \frac{e}{q}, \frac{f}{q}, g,
\frac{h}{q} ) $ in (2.8) which yields the value of $\Phi
(g+, h-)$ in terms of $\Phi_- (g+)$ and $\Phi$. Also,
Lemma 2.1$'$ gives the value of $\Phi(g-, h+)$ in terms
of $\Phi_+(g-)$ and $\Phi$.  Substituting these values
of $\Phi (g+, h-)$ and $\Phi(g-, h+)$ into (2.9) and
simplifying we obtain (2.10).

\bigskip
\noindent{\bf 3. Solutions to a difference equation and a
continued fraction}
\medskip

In the contiguous relation (2.9) of Theorem 2.4, we make 
the replacements
$$
(h,g) \rightarrow \left(hq^{-n}, \frac{s}{h} \qnm \right),
$$
where, to account for the balance condition, we choose
$$
s = {{a^3q^3} \over {bcdef}}.
\leqno{(3.1)}
$$
After renormalization, the above relation can be reduced to the
second order finite difference equation 
$$
\eqalign{
& X_{n+1} -a_n X_n + b_n X_{n-1} =0, \cr
& a_n = A_n +B_n + {{s\qtnm} \over {ah^2}} 
{{(1 - \frac{s}{aq^2})(1-b)(1-c)(1-d)(1-e)(1-f)} \over {
(1- \frac{a}{h} \qnp )(1 - \frac{s}{ah} \qntm)}}, \cr
& b_n = A_{n-1} B_n, \cr
}
\leqno{(3.2)}
$$
where
$$
\eqalign{
& A_n = {{
(1-\frac{s\qnm}{h})
(1-\frac{s\qnm}{ah})
(1-\frac{a\qnp}{bh})
(1-\frac{a\qnp}{ch})
(1-\frac{a\qnp}{dh})
(1-\frac{a\qnp}{eh})
(1-\frac{a\qnp}{fh})
} \over {
(1 - \frac{s\qtn}{h^2})
(1 - \frac{s\qtnm}{h^2})
(1 - \frac{a\qnp}{h})
}} \cr
& B_n = {{
q
(1-\frac{\qn}{h})
(1-\frac{a\qn}{h})
(1-\frac{bs}{ah} \qntm)
(1-\frac{cs\qntm}{ah})
(1-\frac{ds\qntm}{ah})
(1-\frac{es\qntm}{ah})
(1-\frac{fs\qntm}{ah})
} \over {
(1 - \frac{s}{h^2} \qtnm)
(1 - \frac{s\qtnmt}{h^2})
(1 - \frac{s\qnmt}{ah})
}}. \cr
}
$$
It follows that one of the solutions of the finite difference 
equation (3.2) is
$$
\eqalign{
X_n^{(1)} & = 
{{ \left( \frac{s\qtnm}{h^2}, \frac{a\qnp}{h}
\right)_\infty
} \over {
(
\frac{s\qnm}{h},
\frac{s\qnm}{ah},
\frac{a\qnp}{bh},
\frac{a\qnp}{ch},
\frac{a\qnp}{dh},
\frac{a\qnp}{eh},
\frac{a\qnp}{fh}
)_\infty
}} \cr
&  \qquad \times \Phi \left(a;b,c,d,e,f, \frac{s\qnm}{h}, hq^{-n};
q
\right). \cr
}
\leqno{(3.3)}
$$
In fact, (3.3) gives not one solution but seven different
 pairwise linearly independent solutions depending on
the choice of the distinguished parameter out of the
seven parameters $b$, $c$, $d$, $e$, $f$, $s\qnm \slash
h$, $h\qmn$ which define $\Phi$ (see (1.4)). 
In order to distinguish between these
seven solutions we shall write $X_n^{(1),p}$ instead
of $X_n^{(1)}$, $p$ being the distinguished parameter.
Thus $X_n^{(1)}$ represents a set of seven solutions to
equation (3.2).

A second set of seven pairwise linearly independent solutions to
(3.2) is obtained by applying what we call a `reflection
transformation' to (3.2) and (3.3) [9], [10]. That is, in
(3.2) we make the parameter replacements
$$
(a,b,c,d,e,f,s\qnm \slash h, h q^{-n} ) \rightarrow (
q \slash a,
q \slash b,
q \slash c,
q \slash d,
q \slash e,
q \slash f,
hq^{-n+2} \slash s, \qnp \slash h).
\leqno{(3.4)}
$$
It can then be seen that [25]
$$
a_n \rightarrow a_nh^2q^{-2n+1} \slash s, \qquad b_n \rightarrow
b_{n+1} h^4 q^{-4n} \slash s^2.
\leqno{(3.5)}
$$
It is easy to verify that $b_n \rightarrow b_{n+1} h^4 q^{-4n}
\slash s ^2$. In order to 
check $a_n \rightarrow a_nh^2q^{-2n+1} \slash s$,  we need to
verify a polynomial identity of 
degree 14 which we have done on the computer employing MAPLE
software.
 
Having made the above parameter 
replacements we can renormalize so as to arrive back at the
equation (3.2). Consequently, a second 
set of seven solutions is given by
$$
\eqalign{
X_n^{(2)} & = 
{{\left(
\frac{s\qtnm}{h^2}, \frac{sq^n}{ah} \right)_\infty
} \over { \left(
\frac{\qnp}{h},
\frac{aq^n}{h},
\frac{bs\qnm}{ah},
\frac{cs\qnm}{ah},
\frac{ds\qnm}{ah},
\frac{es\qnm}{ah},
\frac{fs\qnm}{ah}
\right)_\infty }} 
\cr
&  \qquad \times \Phi (
q\slash a;
q\slash b,
q\slash c,
q\slash d,
q\slash e,
q\slash f,
h\qmnpt \slash s, \qnp \slash h; q).
}
\leqno{(3.6)}
$$
The seven solutions represented by the different $\Phi$'s of (3.6) 
will be denoted by $X_n^{(2),p}$ 
where $p$ is chosen as the 
distinguished parameter out of the seven parameters 
$q \slash b$,
$q \slash c$,
$q \slash d$,
$q \slash e$,
$q \slash f$, 
$h\qmnpt \slash s$,
$\qnp \slash h$.

We now look
for additional solutions to (3.2) which may be obtained by suitable
parameter replacements. One solution which can be obtained with the
help of the three-term contiguous relation (2.10) derived in 
Theorem~2.5 now follows .

First we interchange $g \leftrightarrow b$ in (2.10) and then make
the
substitutions
$$
\eqalign{
& a=\frac{S^2\qtnmt}{AH^2},\quad
b=\frac{S}{Aq},\quad
c=\frac{BS\qnm}{AH},\quad
d=\frac{CS\qnm}{AH},\cr
& e=\frac{DS}{AH}\qnm,\quad
f=\frac{ES}{AH}\qnm,\quad
g=\frac{FS\qnm}{AH},\quad
h=\frac{S\qnm}{H},\cr
}
\leqno{(3.7)}
$$
ensuring that $bcdefgh=a^3q^2$ and $S=\frac{A^3q^3}{BCDEF}$.

Next we renormalize and get back to equation (3.2) with lower case
letters $a$, $b$, $c$, $d$, $e$, $f$, $s$ replaced by capitals. 
Thus we arrive at a third set of solutions
$$
\eqalign{
X_n^{(3)} & = \left(\frac{s}{aq}\right)^n {{(
\frac{s\qn}{ah},
\frac{s\qtn}{h^2},
\frac{s\qtnm}{h^2} )_\infty
} \over {(
\frac{s\qnm}{h},
\frac{\qnp}{h},
\frac{s^2\qtnm}{ah^2} )_\infty}} \cr
& \quad  \times {{\left(
\frac{s\qn}{bh},
\frac{s\qn}{ch},
\frac{s\qn}{dh},
\frac{s\qn}{eh},
\frac{s\qn}{fh} \right)_\infty
} \over {
\left(
\frac{a\qnp}{bh},
\frac{a\qnp}{ch},
\frac{a\qnp}{dh},
\frac{a\qnp}{eh},
\frac{a\qnp}{fh} \right)_\infty
\left(
\frac{bs\qnm}{ah},
\frac{cs\qnm}{ah},
\frac{ds\qnm}{ah},
\frac{es\qnm}{ah},
\frac{fs\qnm}{ah} \right)_\infty }} \cr
& \quad  \times \Phi \left( \frac{s^2\qtnmt}{ah^2}; \frac{s}{aq}, 
\frac{bs\qnm}{ah},
\frac{cs\qnm}{ah},
\frac{ds\qnm}{ah},
\frac{es\qnm}{ah},
\frac{fs\qnm}{ah},
\frac{s\qnm}{h}; q \right) . \cr
}
\leqno{(3.8)}
$$

We might expect that (3.8) represents seven new solutions 
$X_n^{(3),p}$ where $p$ is chosen as the distinguished parameter 
out of the parameters
$\frac{s}{aq}$,
$\frac{bs\qnm}{ah}$,
$\frac{cs\qnm}{ah}$,
$\frac{ds\qnm}{ah}$,
$\frac{es\qnm}{ah}$,
$\frac{fs\qnm}{ah}$ and
$\frac{s\qnm}{h}$.
However, we find that one of these seven solutions viz. 
$X_n^{(3),\frac{s}{h} \qnm}$ 
is simply a constant multiple of a solution of the first set viz., 
$X_n^{(1), \frac{s}{h}\qnm}$. The actual relation between the 
two solutions is 
$$
X_n^{(1),\frac{s}{h} \qnm} = {{(aq,b,c,d,e,f,h,q \slash h)_\infty}
\over {(
aq\slash b,
aq\slash c,
aq\slash d,
aq\slash e,
aq\slash f,
s\slash aq,
ahq\slash s,
s\slash ah )_\infty }}
X_n^{(3),\frac{s}{h} \qnm} .
\leqno{(3.9)}
$$
Thus (3.8) gives only six new pairwise linearly independent
solutions 
$X_n^{(3),p}$ where $p$ is the 
distinguished parameter chosen out of 
$\frac{s}{aq}$,
$\frac{bs\qnm}{ah}$,
$\frac{cs\qnm}{ah}$,
$\frac{ds\qnm}{ah}$,
$\frac{es\qnm}{ah}$,
$\frac{fs\qnm}{ah}$.

We now apply the reflection transformation (3.4) to the solution 
$X_n^{(3)}$ which enables us to arrive at the solution
$$
\eqalign{
& X_n^{(4)} = \left( \frac{aq^2}{s} \right)^n {{\left(
\frac{s^2\qtnmr}{ah^2} \right)_\infty
} \over {
\left( \frac{s}{ah}\qnm,
\frac{s\qnm}{bh},
\frac{s\qnm}{ch},
\frac{s\qnm}{dh},
\frac{s\qnm}{eh},
\frac{s\qnm}{fh} \right)_\infty}} \cr
& \times \Phi \left(
\frac{ah^2}{s^2}\qtmnpr; 
\frac{aq^2}{s},
\frac{ah\qmnpt}{bs},
\frac{ah\qmnpt}{cs},
\frac{ah\qmnpt}{ds},
\frac{ah\qmnpt}{es},
\frac{ah\qmnpt}{fs},
\frac{h\qmnpt}{s}; q
\right). \cr
}
\leqno{(3.10)}
$$
This gives a set of six new solutions $X_n^{(4),p}$ where $p$ is
chosen out 
of the six parameters 
$\frac{aq^2}{s}$,
$\frac{ah\qmnpt}{bs}$,
$\frac{ah\qmnpt}{cs}$,
$\frac{ah\qmnpt}{ds}$,
$\frac{ah\qmnpt}{es}$,
$\frac{ah\qmnpt}{fs}$. The remaining parameter does not give a new
solution. It is easily seen that 
$X_n^{(4), \frac{h}{s} \qmnpt}$ is a constant multiple of the
already obtained solution 
$X_n^{(2), \frac{h}{s}\qmnpt}$.

A fifth set of solutions to (3.2) can be obtained as follows. After
interchanging 
$g \leftrightarrow b$ 
in the contiguous relation (2.10) of Theorem 2.5, we make the
substitutions
$$
\eqalign{
& 
a = \frac{A\qtnp}{H^2}, \qquad
b = \frac{Aq^2}{S}, \qquad
c = \frac{A\qnp}{BH}, \qquad
d = \frac{A\qnp}{CH}, \qquad \cr
& 
e = \frac{A\qnp}{DH}, \qquad
f = \frac{A\qnp}{EH}, \qquad
g = \frac{A\qnp}{FH}, \qquad
h = \frac{\qnp}{H}, \qquad \cr
}
\leqno{(3.11)}
$$
which again ensures that $bcdefgh=a^3q^2$ and
$S=\frac{A^3q^3}{BCDEF}$. Renormalizing
we obtain (3.2) with lower case letters $a$, $b$, $c$, $d$, $e$,
$f$, $s$ replaced 
by capital letters. This gives another set of solutions 
$$
\eqalign{
& X_n^{(5)} = \left( \frac{aq^2}{s} \right)^n
{{ 
\left( \frac{s}{h^2} \qtnm \right)_\infty
\left( \frac{s}{h^2} \qtn \right)_\infty
\left( \frac{a\qnp}{h} \right)_\infty
} \over {
\left( \frac{a}{h^2} \qtnpt \right)_\infty
\left( \frac{s}{h} \qnm \right)_\infty
\left( \frac{\qnp}{h} \right)_\infty
}} \cr
& \times {{(
\frac{b\qnp}{h},
\frac{c\qnp}{h},
\frac{d\qnp}{h},
\frac{e\qnp}{h},
\frac{f\qnp}{h} )_\infty
} \over { (
\frac{bs}{ah}\qnm,
\frac{cs}{ah}\qnm,
\frac{ds}{ah}\qnm,
\frac{es}{ah}\qnm,
\frac{fs}{ah}\qnm )_\infty (
\frac{a\qnp}{bh},
\frac{a\qnp}{ch},
\frac{a\qnp}{dh},
\frac{a\qnp}{eh},
\frac{a\qnp}{fh} )_\infty }} \cr
& \times \Phi \left( \frac{a\qtnp}{h^2};
\frac{aq^2}{s},
\frac{a\qnp}{bh},
\frac{a\qnp}{ch},
\frac{a\qnp}{dh},
\frac{a\qnp}{eh},
\frac{a\qnp}{fh},
\frac{\qnp}{h};
q \right). \cr
}
\leqno{(3.12)}
$$

Depending on the choice of the distinguished parameter, there are
seven 
solutions given by (3.12). However two of these seven solutions
viz., 
$X_n^{(5), \frac{\qnp}{h}}$ 
and
$X_n^{(5), \frac{aq^2}{s}}$ 
are constant multiples of 
$X_n^{(2), \frac{\qnp}{h}}$ 
and 
$X_n^{(4), \frac{aq^2}{s}}$ 
obtained before. Hence (3.12) gives five new solutions. 

Next, we apply the reflection transformation to the 
solution $X_n^{(5)}$. This leads to the solution 
$$
\eqalign{
& X_n^{(6)} = \left(\frac{s}{aq} \right)^n {{
\left(\frac{a}{h^2}\qtn \right)_\infty} \over {
\left(\frac{a}{h} \qn \right)_\infty \left(
\frac{b}{h}\qn,
\frac{c}{h}\qn,
\frac{d}{h}\qn,
\frac{e}{h}\qn,
\frac{f}{h}\qn \right)_\infty }}\cr
& \times \Phi \left( \frac{h^2}{a}\qmtn; \frac{s}{aq},
\frac{bh}{a}\qmn,
\frac{ch}{a}\qmn,
\frac{dh}{a}\qmn,
\frac{eh}{a}\qmn,
\frac{fh}{a}\qmn, h\qmn; q \right). \cr
}
\leqno{(3.13)}
$$
Again, (3.13) represents only five new solutions because two of the
seven 
solutions viz., 
$X_n^{(6), h\qmn}$
and
$X_n^{(6), \frac{s}{aq} }$ are just constant multiples of the
solutions 
$X_n^{(1), h\qmn}$ 
and $X_n^{(3), \frac{s}{aq}}$ respectively. 

Finally, let us make the following parameter substitutions in the 
contiguous relation (2.9) of Theorem 2.4:
$$
\eqalign{
& a =\frac{B^2}{A}, \qquad
b=B, \qquad
c=\frac{BC}{A}, \qquad
d=\frac{BD}{A}, \qquad \cr
& e=\frac{BE}{A}, \qquad
f=\frac{BF}{A}, \qquad
g=\frac{BS\qnm}{AH}, \qquad
h=\frac{BH\qmn}{A}. \cr
}
\leqno{(3.14)}
$$
Simplifying and renormalizing we again arrive at (3.2) provided we
can verify
an identity in polynomials of degree 14. This also we have done
with the help of
MAPLE software. Consequently we obtain the solution
$$
\eqalign{
& X_n^{(7)} = {{ 
\left( \frac{s}{h^2} \qtnm \right)_\infty
\left( \frac{b}{h} \qnp \right)_\infty
} \over {
\left(
\frac{bs}{ah} \qnm,
\frac{s}{bh} \qnm,
\frac{\qnp}{h},
\frac{a\qnp}{ch},
\frac{a\qnp}{dh},
\frac{a\qnp}{eh},
\frac{a\qnp}{fh}
\right)_\infty }}\cr
& \times \Phi \left(
\frac{b^2}{a}; b,
\frac{bc}{a},
\frac{bd}{a},
\frac{be}{a},
\frac{bf}{a},
\frac{bs\qnm}{ah},
\frac{bh\qmn}{a}; q 
\right).  \cr
}
\leqno{(3.15)}
$$
The solutions 
$X_n^{(7), b}$,
$X_n^{(7), \frac{bs}{ah}\qnm}$,
$X_n^{(7), \frac{bh\qmn}{a}}$
are clearly constant multiples of solutions obtained before. Hence 
$X_n^{(7)}$ represents a set of four new solutions. 

Using reflection on $X_n^{(7)}$ we obtain the solution
$$
\eqalign{
& X_n^{\left(8\right)} = {{ \left(\frac{s}{h^2} \qtn,
\frac{s}{bh}\qn\right)_\infty } \over {
\left(
\frac{s}{h}\qnm,
\frac{a\qnp}{bh},
\frac{b\qn}{h}\right)_\infty \left(
\frac{cs}{ah}\qnm,
\frac{ds}{ah}\qnm,
\frac{es}{ah}\qnm,
\frac{fs}{ah}\qnm\right)_\infty }} \cr
& \times \Phi \left(
\frac{aq}{b^2};
\frac{q}{b},
\frac{aq}{bc},
\frac{aq}{bd},
\frac{aq}{be},
\frac{aq}{bf},
\frac{ah\qmnpt}{bs},
\frac{a\qnp}{bh}; q \right) . \cr
}
\leqno{(3.16)}
$$
$X_n^{(8)}$ gives only four additional solutions because three of
the seven 
solutions viz., 
$X_n^{(8), \frac{q}{b}}$,
$X_n^{(8), \frac{ah}{bs}\qmnpt}$ and 
$X_n^{(8), \frac{a\qnp}{bh}}$ are constant multiples of
previously obtained solutions.

Parameter interchanges $b \leftrightarrow (c,d,e,f)$ in
solutions $X_n^{(7)}$ and $X_n^{(8)}$ give us additional
solutions. We find that there are twelve more solutions
obtained in this manner, the rest being constant multiples
of previous solutions.

Combining (3.3), (3.6), (3.8), (3.10), (3.12), (3.13),
(3.15), (3.16) we have fifty-six pairwise linearly
independent solutions to the three-term recurrence
(3.2). Any three of these fifty-six solutions are connected
by a three-term transformation formula. We shall derive
such a connection formula in \S 6.

For the special case $h=1$ and $s=q,q^2,\dots$, the
solutions $X_n^{(1)}$ and $X_n^{(2)}$ were obtained in
[10].  Writing $h=1$ and $s=q,q^2,\dots$ in (3.2), (3.3)
and (3.6) we get respectively [10, (2.6),
(2.9), (2.12)]  but with a different normalization.

\bigskip
\noindent{\bf Large $n$ asymptotics of the solutions and a
continued fraction}
\medskip

Since $\vert q \vert < 1$, we have from (3.2)
$$
\lim\limits_{n \rightarrow \infty} a_n = 1 + q, \qquad
\lim\limits_{n \rightarrow \infty} b_n=q.
\leqno{(3.17)}
$$
It follows that (3.2) has solutions whose large $n$ asymptotics is
either
$$
X_n \mathop{\approx}\limits^{n \rightarrow \infty} {\rm\ const.\ }
\leqno{(3.18a)}
$$
or
$$
X_n \mathop{\approx}\limits^{n \rightarrow \infty} {\rm\ const.\ }
q^n .
\leqno{(3.18b)}
$$
The latter characterizes the minimal solution to (3.2). 
We therefore proceed to examine the asymptotics 
of the solutions obtained above. From (3.3) we have in 
a straightforward manner
$$
\eqalign{
& X_n^{(1), b} \mathop{\approx}\limits^{n \rightarrow \infty} 
{\w} (a;b,c,d,e,f; \frac{s}{aq} ) \cr
& + {{(aq,c,d,e,f, \frac{b}{a},
\frac{bq}{c},
\frac{bq}{d},
\frac{bq}{e},
\frac{bq}{f}, h,
\frac{q}{h},
\frac{bhq^2}{s},
\frac{s}{bhq})_\infty} \over {(
\frac{b^2q}{a},
\frac{bc}{a},
\frac{bd}{a},
\frac{be}{a},
\frac{bf}{a},
\frac{a}{b},
\frac{aq}{c},
\frac{aq}{d},
\frac{aq}{e},
\frac{aq}{f},
\frac{bh}{a},
\frac{aq}{bh},
\frac{ahq^2}{s},
\frac{s}{ahq}
)_\infty }} \cr
& \times {\w} \left(
\frac{b^2}{a}; b,
\frac{bc}{a},
\frac{bd}{a},
\frac{be}{a},
\frac{bf}{a};
\frac{s}{aq} \right) , \qquad \left\vert \frac{s}{aq} \right\vert
<1. \cr
}
\leqno{(3.19)}
$$
We will have analogous results with $c$, $d$, $e$ or $f$ as
distinguished 
parameters. We also easily have 
$$
X_n^{(1), \frac{s}{h} \qnm} \mathop{\approx}\limits^{n \rightarrow
\infty} 
{\w} (a;b,c,d,e,f; \frac{s}{aq}), \qquad \left\vert\frac{s}{aq}
\right\vert <1 .
\leqno{(3.20)}
$$
In order to calculate the asymptotics of $X_n^{(1), h\qmn}$,
we first apply the transformation [7, (2.30), p.~56]
to the $\Phi^{(h\qmn)}$ in (3.3). While working out the
asymptotics of the two resulting   $\tphia$ series, some
care is necessary, since when $n$ is large, the terms
near  two parts of the series are important (see Bailey
[4], p.~128). Thus it is convenient to break each $\tphia$
series into three parts
$$
\sum\limits_{k=0}^{n \slash 2} +
\sum\limits_{k=n \slash 2}^{3 n \slash 2} +
\sum\limits_{k=3 n \slash 2}^{\infty} ,
$$
the third sum tending to $0$ as $ n \rightarrow \infty$. The first 
summation gives a ${_4 \phi_3}$ while the second gives a ${_4
\psi_4}$.
We obtain 
$$
\eqalign{
& X_n^{(1), h\qmn} \mathop{\approx}\limits^{n \rightarrow \infty} 
{{(aq,c,
\frac{aq}{db}, \frac{aq}{eb}, \frac{aq}{fb}, \frac{sb}{aq},
\frac{h}{a},
\frac{hq}{c}, \frac{ahq^2}{bs}, \frac{bdh}{a}, \frac{beh}{a},
\frac{bfh}{a})_\infty } \over {(
\frac{c}{b}, \frac{aq}{b}, \frac{aq}{d}, \frac{aq}{e},
\frac{aq}{f},
\frac{s}{aq}, \frac{bhq}{c}, \frac{bh}{a}, \frac{ahq^2}{s},
\frac{dh}{a},
\frac{eh}{a},
\frac{fh}{a} )_\infty }} \cr
& \times {{(
\frac{aq}{h}, \frac{c}{h}, \frac{bs}{ahq}, \frac{aq}{bdh},
\frac{aq}{beh},
\frac{aq}{bfh} )_\infty } \over {(
\frac{c}{bh}, \frac{aq}{bh}, \frac{s}{ahq}, \frac{aq}{dh},
\frac{aq}{eh},
\frac{aq}{fh} )_\infty }} \left[ {_4 \phi_3} \left( 
{{
\frac{aq}{cd}, \frac{aq}{ce},
\frac{aq}{cf}, b} \atop {
\frac{bq}{c}, \frac{bs}{aq},
\frac{aq}{c} }} ; q \right) \right. \cr
& \left. - \frac{c}{bh} {{(
\frac{cq}{bh}, \frac{q}{h}, \frac{cs}{ahq}, \frac{aq}{bh},
\frac{aq}{dc},
\frac{aq}{ec},
\frac{aq}{fc}, b)_\infty } \over {(
\frac{aq}{bdh}, \frac{aq}{beh}, \frac{aq}{bfh}, \frac{c}{h},
\frac{bq}{c},
\frac{bs}{aq},
\frac{aq}{c}, q)_\infty }} {_4 \psi_4} \left( {{
\frac{bh}{c}, h, \frac{ahq^2}{cs}, \frac{bh}{a} } \atop {
\frac{bdh}{a}, \frac{beh}{a}, \frac{bfh}{a},
\frac{hq}{c} }} ; q \right) \right] \cr
& + {\rm \ idem\ } (b;c) . \cr
}
\leqno{(3.21)}
$$
The asymptotics of the seven solutions represented by $X_n^{(1)}$
is thus completely given by (3.19), (3.20) and (3.21). There are
similar results
for $X_n^{(2)}$ solutions. We easily have 
$$
\eqalign{
& X_n^{(2), \frac{q}{b}} \mathop{\approx}\limits^{n \rightarrow
\infty} 
\left[ {\w} \left( \frac{q}{a};
\frac{q}{b},
\frac{q}{c},
\frac{q}{d},
\frac{q}{e},
\frac{q}{f};
\frac{aq^2}{s} \right) \right. \cr
& + {{ (
\frac{q^2}{a},
\frac{a}{b},
\frac{q}{c},
\frac{q}{d},
\frac{q}{e},
\frac{q}{f},
\frac{cq}{b},
\frac{dq}{b},
\frac{eq}{b},
\frac{fq}{b},
\frac{hq^2}{s},
\frac{s}{hq},
\frac{hq}{b},
\frac{b}{h} )_\infty } \over {(
\frac{aq^2}{b^2},
\frac{b}{a},
\frac{cq}{a},
\frac{dq}{a},
\frac{eq}{a},
\frac{fq}{a},
\frac{aq}{bc},
\frac{aq}{bd},
\frac{aq}{be},
\frac{aq}{bf},
\frac{hq}{a},
\frac{a}{h},
\frac{ahq^2}{bs},
\frac{bs}{ahq} )_\infty }} \cr
& \left.  \times {\w} \left( 
\frac{aq}{b^2};
\frac{q}{b},
\frac{aq}{bc},
\frac{aq}{bd},
\frac{aq}{be},
\frac{aq}{bf};
\frac{aq^2}{s} \right) \right], \qquad 
\left\vert \frac{aq^2}{s} \right\vert <1 , \cr
}
\leqno{(3.22)}
$$
and 
$$
X_n^{(2), \frac{\qnp}{h}} \mathop{\approx}^{n \rightarrow \infty} 
{\w} \left( 
\frac{q}{a};
\frac{q}{b},
\frac{q}{c},
\frac{q}{d},
\frac{q}{e},
\frac{q}{f};
\frac{aq^2}{s} 
\right), \qquad \left\vert 
\frac{aq^2}{s}  \right\vert < 1.
\leqno{(3.23)}
$$
The asymptotics of 
$X_n^{(2), \frac{h}{s} \qmnpt }$ 
is worked out exactly as we have done for 
$X_n^{(1), h \qmn }$. The result being cumbersome is not 
being given here. 

For $X_n^{(3)}$, $X_n^{(4)}$ we have the following results 
which take care of all the related solutions.
$$
\eqalign{
& X_n^{(3), \frac{bs}{ah}\qnm} \mathop{\approx}^{n \rightarrow
\infty} 
{{(
\frac{s}{aq},
\frac{bq}{c},
\frac{bq}{d},
\frac{bq}{e},
\frac{bq}{f},
\frac{bq}{a},
\frac{bhq}{s},
\frac{s}{bh} )_\infty } \over { (
\frac{b^2q}{a},
\frac{bc}{a},
\frac{bd}{a},
\frac{be}{a},
\frac{bf}{a}, b,
\frac{bh}{a},
\frac{aq}{bh} )_\infty }} \cr
& \times {\w} \left(
\frac{b^2}{a};
\frac{bc}{a},
\frac{bd}{a},
\frac{be}{a},
\frac{bf}{a}, b; \frac{s}{aq} \right),
\qquad \left\vert \frac{s}{aq}  \right\vert < 1. \cr
}
\leqno{(3.24)}
$$

$$
\eqalign{
& X_n^{(3), \frac{s}{aq}} \mathop{\approx}^{n \rightarrow \infty} 
 {{ ( 
\frac{a}{b},
\frac{q}{c},
\frac{q}{d},
\frac{q}{e},
\frac{bs}{aq},
\frac{cs}{aq},
\frac{ds}{aq},
\frac{es}{aq},
\frac{h^2q}{s},
\frac{s}{h^2},
\frac{hq}{f}, \frac{f}{h} )_\infty } \over { (
\frac{s}{f},\frac{f}{a}, h,
\frac{q}{h},
\frac{bh}{a},
\frac{aq}{bh},
\frac{ch}{a},
\frac{aq}{ch},
\frac{dh}{a},
\frac{aq}{dh},
\frac{eh}{a},
\frac{aq}{eh} )_\infty }}  \cr
& \times {\w} \left(
\frac{s}{fq};
\frac{s}{aq},
\frac{aq}{bf},
\frac{aq}{cf},
\frac{aq}{df},
\frac{aq}{ef}; f \right) \cr
& + {{(
\frac{aq}{bf},
\frac{aq}{cf},
\frac{aq}{df},
\frac{aq}{ef},
\frac{bfs}{a^2q},
\frac{cfs}{a^2q},
\frac{dfs}{a^2q},
\frac{efs}{a^2q},
\frac{h^2q}{s},
\frac{s}{h^2},
\frac{hq}{a},
\frac{a}{h} )_\infty } \over { (
\frac{fs}{a^2},
\frac{a}{f},
\frac{bh}{a},
\frac{aq}{bh},
\frac{ch}{a},
\frac{aq}{ch},
\frac{dh}{a},
\frac{aq}{dh},
\frac{eh}{a},
\frac{aq}{eh},
\frac{fh}{a},
\frac{aq}{fh} )_\infty }} \cr
& \times {\w} \left( \frac{sf}{a^2q}; \frac{s}{aq},
\frac{q}{b},
\frac{q}{c},
\frac{q}{d},
\frac{q}{e}; f \right), \qquad \vert f \vert < 1. \cr
}
\leqno{(3.25)}
$$
In order to obtain asymptotics of $X_n^{(4), \frac{aq^2}{s}}$, we 
first apply the transformation ([7], (2.30), p.~56) to the 
$\Phi$ in the solution and then work out the asymptotics. We get
$$
\eqalign{
& X_n^{(4), \frac{aq^2}{s}} \mathop{\approx}^{n \rightarrow \infty} 
\cr
& \left[ {{(b,c,d,e,
\frac{aq^3}{bs}, \frac{aq^3}{cs}, \frac{aq^3}{ds}, \frac{aq^3}{es},
\frac{ah^2q^4}{s^2}, \frac{s^2}{ah^2q^3}, \frac{ahq^2}{fs},
\frac{fs}{ahq} )_\infty } \over {(
\frac{a}{f}, \frac{fq^3}{s}, \frac{ahq^2}{s}, \frac{bhq^2}{s},
\frac{chq^2}{s}, \frac{dhq^2}{s}, \frac{ehq^2}{s}, \frac{s}{ahq},
\frac{s}{bhq}, \frac{s}{chq}, \frac{s}{dhq},
\frac{s}{ehq} )_\infty }} \right. \cr
& \times {\w} \left(
\frac{fq^2}{s}; \frac{aq^2}{s}, \frac{bf}{a}, \frac{cf}{a},
\frac{df}{a}, \frac{ef}{a};
\frac{q}{f} \right) \cr
& + {{(
\frac{bf}{a}, \frac{cf}{a}, \frac{df}{a}, \frac{ef}{a},
\frac{a^3q^3}{bfs}, \frac{a^3q^3}{cfs}, \frac{a^3q^3}{dfs}, 
\frac{a^3q^3}{efs}, \frac{ah^2q^4}{s^2}, \frac{s^2}{ah^2q^3},
\frac{hq^2}{s},
\frac{s}{hq} )_\infty } \over { (
\frac{a^2q^3}{fs}, \frac{f}{a}, \frac{hbq^2}{s}, \frac{hcq^2}{s},
\frac{hdq^2}{s}, \frac{heq^2}{s}, \frac{hfq^2}{s}, \frac{s}{hbq},
\frac{s}{hcq}, \frac{s}{hdq}, \frac{s}{heq},
\frac{s}{hfq} )_\infty }} \cr
& \left. \vphantom{{\frac{a^3q^3}{cfs}}\over{\frac{a^3q^3}{dfs}}}
\times {\w} \left(
\frac{a^2q^2}{fs}; \frac{aq^2}{s},b,c,d,e; \frac{q}{f} \right)
\right], \qquad
\left\vert \frac{q}{f} \right\vert < 1. \cr
}
\leqno{(3.26)}
$$
The calculation of asymptotics of $X_n^{(4), \frac{ah\qmnpt}{bs}}$
is similar to that of $X_n^{(1), h\qmn}$. 

Coming to solutions $X_n^{(5)}$ and $X_n^{(6)}$ we find that 
$$
\eqalign{
& X_n^{(5), \frac{a\qnp}{bh}} \mathop{\approx}^{n \rightarrow
\infty}
{{( \frac{aq^2}{s},
\frac{cq}{b},
\frac{dq}{b},
\frac{eq}{b},
\frac{fq}{b},
\frac{aq}{b},
\frac{h}{b},
\frac{bq}{h} )_\infty } \over { (
\frac{aq^2}{b^2},
\frac{aq}{bc},
\frac{aq}{bd},
\frac{aq}{be},
\frac{aq}{bf},
\frac{q}{b},
\frac{ahq^2}{bs},
\frac{bs}{ahq} )_\infty }} \cr
& \times {\w} \left( 
\frac{aq}{b^2};
\frac{aq}{bc},
\frac{aq}{bd},
\frac{aq}{be},
\frac{aq}{bf},
\frac{q}{b};
\frac{aq^2}{s} \right), \qquad \left\vert \frac{aq^2}{s}
\right\vert < 1, \cr
}
\leqno{(3.27)}
$$
while $X_n^{(6), \frac{bh}{a}\qmn}$ is worked out as $X_n^{(1),
h\qmn}$
and produces a similar expression giving constant asymptotics. 
The remaining solutions are structurally similar to the 
solutions for which we have calculated the asymptotics above.
Thus we find that all the fifty-six solutions yield constant
asymptotics
as $n \rightarrow \infty$.

In order to obtain a minimal solution to (3.2) we make use of 
(3.20) and (3.23). We write 
$$
W_1 \coloneq {\w} \left(a;b,c,d,e,f; \frac{s}{aq} \right),
\qquad \left\vert \frac{s}{aq} \right\vert <1
\leqno{(3.28)}
$$
and its analytic continuation otherwise and 
$$
W_2 \coloneq {\w} \left(\frac{q}{a};
\frac{q}{b},\frac{q}{c},\frac{q}{d},\frac{q}{e}, 
\frac{q}{f}; \frac{aq^2}{s} \right),
\qquad \left\vert \frac{aq^2}{s} \right\vert <1
\leqno{(3.29)}
$$
and its analytic continuation otherwise. Define
$$
\eqalign{
X_n^{(\min)} & \coloneq W_2X_n^{(1), \frac{s}{h}\qnm} - 
W_1 X_n^{(2), \frac{\qnp}{h}} \cr
& \mathop{\approx}^{n \rightarrow \infty} {\rm \ const.\ } q^n 
\qquad {\rm from\ (3.20)\ and\ (3.23).}
}
\leqno{(3.30)}
$$
The actual value of the constant in (3.30) will be given 
in section 7. It follows from  (3.30) and (3.20) that 
$$
\lim\limits_{n \rightarrow \infty} 
{{X_n^{(\min)} } \over {X_n^{(1), \frac{s}{h} \qnm}}} =0, \qquad
\left\vert \frac{s}{q} \right\vert < \vert a \vert < \left\vert
\frac{s}{q^2} \right\vert.
\leqno{(3.31)}
$$
Thus $X_n^{(\min)}$ given by (3.30) is a minimal solution of (3.2).
An application of Pincherle's theorem [8], [19] then leads to the
following
result.

\medskip
\noindent{\bf Theorem 3.1.}
{\it For $a_n$, $b_n$ defined by {\rm (3.2)} the following 
continued fraction representation holds true:\/}
$$
\eqalign{
& \frac{1}{a_0} \ \lowminus
\frac{b_1}{a_1} \ \lowminus
\frac{b_2}{a_2} \ \lowminus \ \lowdots  
 = {{W_2X_0^{(1), \frac{s}{hq}} - W_1X_0^{(2), \frac{q}{h}}  
} \over { b_0 \left(
W_2X_{-1}^{(1), \frac{s}{hq^2}} - W_1X_{-1}^{(2), \frac{1}{h}}  
\right) }} \cr
& = \frac{(1-\frac{s}{h^2q} )(1- \frac{s}{h^2q^2})}{q} \cr
& \times \left[ W_2 {{(\frac{aq}{h})_\infty } \over { (
\frac{s}{hq}, \frac{s}{ahq}, \frac{aq}{bh}, \frac{aq}{ch},
\frac{aq}{dh}, \frac{aq}{eh},
\frac{aq}{fh})_\infty }} 
\Phi^{(\frac{s}{hq})} \left(a;b,c,d,e,f, \frac{s}{hq},h;q
\right)  \right. \cr
&\left.  -W_1 {{(\frac{s}{ah})_\infty } \over { (
\frac{q}{h}, \frac{a}{h}, \frac{bs}{ahq}, \frac{cs}{ahq},
\frac{ds}{ahq}, \frac{es}{ahq},
\frac{fs}{ahq} )_\infty }} 
\Phi^{(\frac{q}{h})} \left(
\frac{q}{a}; \frac{q}{b}, \frac{q}{c}, \frac{q}{d}, \frac{q}{e},
\frac{q}{f}, \frac{hq^2}{s},
\frac{q}{h}; q \right) 
\vphantom{{{\frac{q}{e}}\over{\frac{q}{f}}}} \right]
\cr
&\bigslash  \left[ W_2 {{
(1 -\frac{1}{h}) (1 -\frac{bs}{ahq^2}) (1 -\frac{cs}{ahq^2})
(1 -\frac{ds}{ahq^2}) (1 -\frac{es}{ahq^2}) (1 -\frac{fs}{ahq^2})
(\frac{a}{h})_\infty } \over {(
\frac{s}{hq}, \frac{s}{ahq^2}, \frac{aq}{bh}, \frac{aq}{ch},
\frac{aq}{dh}, \frac{aq}{eh},
\frac{aq}{fh} )_\infty }} \right. \cr
& \times \Phi^{(\frac{s}{hq^2})} \left(a;b,c,d,e,f, \frac{s}{hq^2},
hq; q \right) \cr
& - W_1 {{
(1- \frac{s}{hq^2}) (1- \frac{a}{bh}) (1- \frac{a}{ch}) (1-
\frac{a}{dh})
(1- \frac{a}{eh}) (1- \frac{a}{fh})
(\frac{s}{ahq})_\infty } \over { (
\frac{q}{h}, \frac{a}{qh}, \frac{bs}{aqh}, \frac{cs}{aqh},
\frac{ds}{aqh},
\frac{es}{aqh},
\frac{fs}{aqh} )_\infty }} \cr 
& \left. \times \Phi^{(\frac{1}{h})} \left(
\frac{q}{a}; \frac{q}{b}, \frac{q}{c}, \frac{q}{d}, \frac{q}{e},
\frac{q}{f}, \frac{hq^3}{s},
\frac{1}{h};
q \right) 
\vphantom{{{\frac{q}{e}}\over{\frac{q}{f}}}}\right] .
\cr
}
\leqno{(3.32)}
$$
For the special case $h=1$ we have 

\medskip
\noindent{\bf Corollary 3.2.} {\it If $a_n$, $b_n$ are
defined by (3.2) with $h=1$, then the following continued
fraction representation holds true\/} (Masson [25]):
$$
\eqalign{
& {1 \over {a_0}}  \lowminus {{b_1} \over {a_1}}
\lowminus {{b_2} \over {a_2}} \lowminus \lowdots \cr
& = {{
\left( 1 - \frac{s}{q} \right)
\left( 1 - \frac{a}{q} \right)
} \over { q
\left( 1 - \frac{s}{aq} \right)
\left( 1 - \frac{a}{b} \right)
\left( 1 - \frac{a}{c} \right)
\cdots 
\left( 1 - \frac{a}{f} \right) }} \cr
& \qquad \times \left[ {\vphantom{ {{\left( \frac{f}{q}
\right)_\infty} \over {\left(\frac{q}{f}\right)_\infty}}
}}\W \left( \frac{q}{a} ;
q, \frac{q^2}{s}, \frac{q}{b}, \cdots , \frac{q}{f} ;
q \right) +\Pi_1 \, \W \left( aq; q, \frac{aq^2}{s},
\frac{aq}{b} , \cdots , \frac{aq}{f} ; q \right) \right. \cr
& \qquad\qquad \left. - {{
\left( q,a,aq, \frac{bs}{aq} , \cdots , \frac{fs}{aq}
\right)_\infty} \over { \left( \frac{s}{q},
\frac{s}{a}, \frac{s}{aq}, \frac{aq}{b}, \cdots ,
\frac{aq}{f} \right)_\infty }} {{W_2} \over
{W_1}}\right] \bigslash (1+\Pi_0) \cr
\Pi_n & \coloneq {{\left( \frac{q^2}{a}, \frac{q}{b},
\cdots , \frac{q}{f}, \frac{q^{3-n}}{s}, a\qnmo,s\qtnmt
, b\qn , \dots  f\qn \right)_\infty } \over { \left(
a\qtn , \frac{\qomn}{a}, \frac{bq}{a}, \dots 
\frac{fq}{a}, \frac{aq^2}{s}, \frac{s\qnmo}{a},
\frac{a\qn}{b}, \dots \frac{a\qn}{f} \right)_\infty
}} . \cr
}
\leqno{(3.33)}
$$

\medskip
\noindent{\bf Proof:} We take limit of (3.32) as $h
\rightarrow 1$. 
A special case of Corollary 3.2 for $s = q^3, q^4,
\dots$ was obtained in [10, (5.1), p.~438]. 

In the terminating case of Theorem 3.1 we have

\medskip
\noindent{\bf Corollary 3.3.} {\it If $h=1$  and one of 
$\frac{aq}{b}, \frac{aq}{c}, \frac{aq}{d}, \frac{aq}{e}, 
\frac{aq}{f} = \qmN$, $N=0,1,\dots$, then \/}
$$
\eqalign{
& {1 \over {a_0}} \lowminus {{b_1} \over {a_1}} \lowminus
{{b_2} \over {a_2}} \lowminus \lowdots \lowminus
{{b_N} \over {a_N}}  \cr
& = \frac{aq}{s} {{\left( 1 - aq \right) } \over { 
(1-b)
(1-c)
(1-d)
(1-e)
(1-f) }} \W \left( aq; q,
\frac{aq^2}{s},\frac{aq}{b},\dots \frac{aq}{f} ; q
\right) . \cr
}
\leqno{(3.34)}
$$

\medskip \noindent{\bf Proof:} We take the limit of (3.33)
as, say, $\frac{aq}{f} \rightarrow \qmN$ (see Masson [25]).

It is this last Corollary which yields the rational
biorthogonality discussed in Section~4.

The case $s=q^2$ is Watson's $q$-analogue of Ramanujan's Entry 40
[10].

\bigskip
\noindent{\bf 4. Discrete Rational Biorthogonality}
\medskip

The continued fraction (3.34) can be put into the form of
an $R_{II}$-fraction [17]. Associated with this
$R_{II}$-fraction there is an explicit system of
biorthogonal rational functions [29], [26], [25]. In this
section we derive such a biorthogonal system using the
methods in [17]. This gives the top level of an Askey
type scheme of discrete rational biorthogonality for
which we also outline five different next level limits
and a $q \rightarrow 1$ limit. For the Askey-scheme of
hypergeometric orthogonal polynomials see [20].
\medskip
{\bf 4 a). The top $\tphia$ level.}
In the three term recurrence (3.2) we make the replacements

$$
b \rightarrow be^{-\xi}, \qquad
c \rightarrow \mu be^\xi , \qquad h \rightarrow 1
\leqno{(4.1)}
$$
with
$$
x:= (e^\xi + e^{-\xi} \slash \mu ) \slash 2 .
\leqno{(4.2)}
$$
The coefficients $a_n$ and $b_n$ in (3.2) then become
linear and quadratic functions of $x$ respectively and
the recurrence takes the form
$$
X_\npo (x) - (u_\npo x + v_\npo)X_n (x) + \gamma_n (x -
\alpha_\npo)(x - \beta_\npo)X_n (x) = 0 .
\leqno{(4.3)}
$$
After a renormalization this becomes the recurrence 
$$
Y_\npo (x) - (x-c_\npo)Y_n (x) +\lambda_\npo
(x-\alpha_\npo)(x-\beta_\npo) Y_\nmo (x)=0
\leqno{(4.3')}
$$
associated with the $R_{II}$-fraction
$$
R_{II} (x) = {1 \over {x-c_1}} \lowminus {{\lambda_2
(x-\alpha_2)(x -\beta_2)} \over {x-c_2}} \lowminus
{{\lambda_3(x-\alpha_3)(x-\beta_3)} \over {x-c_3}}
\lowminus \lowdots .
\leqno{(4.4)}
$$
The connection between (4.3) and (4.3$'$) is given by 
$$
c_\npo = - {{v_\npo} \over {u_\npo}}, \qquad \lambda_\npo
= {{\gamma_n} \over {u_\npo u_n}} , \qquad Y_n = {{X_n}
\over {\prod_{k=-1}^n u_k}} ,
\leqno{(4.5)}
$$
but the renormalization factor $\prod_{k=-1}^n u_k$ will
not enter our final formulas. We may explicitly calculate
$c_\npo$, $\lambda_\npo$, $\alpha_\npo$ and $\beta_\npo$
from the original $a_n$ and $b_n$ but for our purposes we
only note that
$$
\eqalign{
\gamma_n & = -4\mu s\qtnmt (1-\qn)(1-s\qnmt)
\left( 1 - {a \over d} \qn\right) 
\left( 1 - {a \over e} \qn\right) 
\left( 1 - {a \over f} \qn\right)  \cr
& \qquad \times 
\left(1 - {{ds} \over a} \qnmt\right)
\left(1 - {{es} \over a} \qnmt\right)
\left(1 - {{fs} \over a} \qnmt\right) , \cr
s & = {{a^3 q^3} \over {\mu b^2 def}}, \cr
}
\leqno{(4.6)}
$$
and the `interpolation points' $\alpha_\npo$,
$\beta_\npo$ are given by 
$$
\eqalign{
\alpha_\npo & = {\left({{bs} \over a} \qnmt + {a \over
{bs\mu}} \qtmn\right)} \slash 2 , \cr
\beta_\npo & = {{\left({{a} \over {b\mu}} \qn + {b \over
{a}} \qmn\right)} \slash 2} . \cr
}
\leqno{(4.7)}
$$

In order to have explicit singularities in the continued
fraction (4.4), we will need termination. For this we now
choose
$$
f = a\qNpo, \qquad N=0,1,\dots
\leqno{(4.8)}
$$
so that $\lambda_{N+2} =0$. We can then use the result
of Corollary 3.4 with the replacement (4.1) to obtain
$$
\eqalign{
R_{II} (x) & = {{u_0 \mu b^2 de\qN(1-aq)} \over {a
(1-be^{-\xi}) (1-\mu be^{\xi}) (1-d)(1-e)(1-a\qNpo)}} \cr
& \qquad \times  \W \left( aq; q, {{aqe^\xi} \over {b}},
{{aqe^{-\xi}} \over {\mu b}}, {{aq} \over {d}}, {{aq} \over {e}},
{{\mu b^2de\qN} \over {a}}, \qmN ;q \right) \cr
& = \sum_{k=0}^N {{R_k} \over {x-x_k}}, \qquad x =
{{\left(e^\xi + \mu^{-1} e^{-\xi} \right) } \over 2} 
. \cr
}
\leqno{(4.9)}
$$
Now an explicit polynomial solution to (4.3) is given by
(3.3) with the replacement (4.1) and identification (4.2).
Consequently the explicit polynomial solution of the
first kind to (4.3$'$) is given by 
$$
P_n (x) = {{\left( 
{{s} \over {q}},
{{s} \over {aq}},
{{aqe^{\xi}} \over {b}},
{{aqe^{-\xi}} \over {\mu b}},
{{aq} \over {d}},
{{aq} \over {e}}, \qmN \right)_n U_n (x)} \over
{(aq)_n\left( {s \over q} \right)_{2n} \prod_{k = -1}^n
u_k }} 
\leqno{(4.10)}
$$
where 
$$
\eqalign{
U_n (x) & = \W \left( a; be^{-\xi}, \mu be^{\xi},
d,e,a\qNpo , s\qnmo, \qmn ; q \right), \cr
s & = {{a^2 \qtmN} \over {\mu b^2 de}} . \cr
}
\leqno{(4.11)}
$$
The polynomial $P_n (x)$ and hence the rational function
$U_n (x)$ satisfy a finite discrete orthogonality
associated with the poles and residues of the continued
fraction (4.9). From the results in Ismail and Masson [17], 
we know this to be 
$$
\sum_{k=0}^N {{P_n (x_k)x_k^m R_k} \over {\prod_{j=2}^n
(x_k -\alpha_{j+1} )(x_k - \beta_{j+1})}} = 0 , \qquad 0
\leq m < n
\leqno{(4.12)}
$$
or, equivalently, 
$$
\sum_{k=0}^N {{U_n (x_k)x_k^m R_k} \over {\prod_{j=1}^n
(x_k -\alpha_{j+1} )}} = 0 , \qquad 0 \leq m < n .
\leqno{(4.12')}
$$
It remains to calculate the $x_k$ and $R_k$ and then
adjust the partial orthogonality (4.12$'$) to obtain an explicit
rational biorthogonality.

The singular points $x_k$ are easily calculated from
(4.9). They are given by the zeros of the denominator
factors $(be^{-\xi},\mu be^\xi)_{k+1}$, $k=0,1,\dots ,N$.
This gives poles at $e^\xi = b\qk$ or ${{\qmk} \slash {\mu
b}}$, $k=0,1,\dots ,N$. That is, 
$$
x_k = {{\left( b\qk + {{\qmk} \slash {b\mu}} \right)} \over
2} , \qquad k=0,1,\dots ,N .
\leqno{(4.13)}
$$

The explicit calculation of the residues $R_k$ is
possible because they are expressed in terms of a
terminating very well poised $\ephis$ for which one may
use the Jackson summation formula. In detail we have 
$$
\eqalign{
R_k & = C\qmk \lim_{e^\xi \rightarrow b\qk} {{(1-\qk
e^{-\xi} b)(1-\mu q^k e^\xi b)} \over {(1-e^{-\xi}
b)(1 - \mu e^\xi b)}} \cr
& \qquad \times \W \left( aq; q, 
{{aqe^{\xi}} \over {b}},
{{aqe^{-\xi}} \over {\mu b}},
{{aq} \over {d}},
{{aq} \over {e}},
{{\mu b^2 de\qN} \over {a}},
\qmN ; q \right) , \cr
C & = - {{u_0 b de (1-aq) \qN} \over
{(1-d)(1-e)(1-a\qNpo)}} . \cr
}
$$
This limit gives 
$$
\eqalign{
R_k & =  C (-1)^k q^{ {{k(k-1)} \slash 2} } {{(1-\mu \qtk
b^2)(1-a\qtkpo) } \over {(1-\mu \qk b^2)(1-aq)}} \cr
& \qquad \times {{\left( a\qkpo , \frac{a}{\mu b^2} \qmkp ,
\frac{aq}{d}, \frac{aq}{e}, \qmN , \frac{\mu b^2 de}{a}
\qN\right)_k} \over { \left( \mu b^2 \qkpo , dq, eq,
a\qNpt , \frac{a^2 \qtmN}{\mu b^2 de} ,q \right)_k }} \cr
& \qquad \times  \w \left( a\qtkpo ; \frac{aq}{\mu b^2} ,
\frac{a\qkpo}{d}, \frac{a\qkpo}{e} , \frac{\mu b^2 de
\qNpk}{a} , \qmNpk ; q \right) . \cr
}
$$
The above $\w$ can be summed using the Jackson
summation [7, ($\caprom 2$.~22), p.~238] 
to give, after some simplification,
$$
\eqalign{
R_k & = C {{\left( aq^2, 
\frac{\mu b^2 d}{a}, 
\frac{\mu b^2 e}{a}, 
\frac{de}{a} \right)_N} \over {\left( \mu b^2 q, dq, eq,
\frac{\mu b^2 de}{a^2} \right)_N}} \cr
& \qquad \times {{\left( \mu b^2 ,
q\sqrt{\mu b^2},
- q\sqrt{\mu b^2},
\frac{\mu b^2}{a},
\frac{aq}{d},
\frac{aq}{e},
\frac{\mu b^2 de \qN}{a},
\qmN \right)_k } \over { \left( 
\sqrt{\mu b^2},
- \sqrt{\mu b^2}, aq, 
\frac{\mu b^2 d}{a},
\frac{\mu b^2 e}{a},
\frac{a}{de} \qomN, \mu b^2 \qNpo , q \right)_k }} . \cr
}
\leqno{(4.14)}
$$
The partial orthogonality (4.12$'$) is now explicit. We
next modify it to obtain an explicit biorthogonality.

Note that in the orthogonality expression (4.12$'$) there
is a denominator factor
$$
x_k -\alpha_2 = \frac{\qmk}{2b\mu} \left(1 - \frac{a}{de}
q^{1-N+k} \right) \left( 1 - \frac{\mu b^2 de}{a}
q^{k+N-1} \right) .
\leqno{(4.15)}
$$
We incorporate this factor into $R_k$ and define a new
weight 
$$
\eqalign{
\omega_k & = {{R_k\left( 1- \frac{a}{de} \qmNpo \right)
\left( 1 - \mu\frac{b^2de}{a} \qNmo \right) } \over { R_0
(x_k -\alpha_2) 2 \mu b}} \cr
& =  {{\left(\mu b^2, 
q \sqrt{\mu b^2},
- q \sqrt{\mu b^2},
\frac{\mu b^2}{a},
\frac{aq}{d},
\frac{aq}{e},
\frac{\mu b^2 de \qNmo}{a} \right)_k \qk } \over {
\left(  \sqrt{\mu b^2},
- \sqrt{\mu b^2}, aq,
\frac{\mu b^2 d}{a},
\frac{\mu b^2 e}{a},
\mu b^2 \qNpo, \frac{a}{de} \qmNpt , q \right)_k }} . \cr
}
\leqno{(4.16)}
$$
The orthogonality (4.12$'$) can then be restated as
$$
\sum_{k=0}^N {{U_n (x_k) Q_m (x_k) \omega_k } \over
{\prod_{j=2}^N (x_k -\alpha_{j+1} ) }} = 0, \qquad 0 \leq
m <n 
\leqno{(4.17)}
$$
where $Q_m (x)$ is any polynomial of degree $m$.
We will now use a symmetry of $w_k$ to obtain a full
biorthogonality.

Consider the parameter interchange
$$
a \leftrightarrow {{a\qomN} \over {de}} 
\leqno{(4.18)}
$$
with $\mu$, $b$, $\frac{aq}{d}$, $\frac{aq}{e}$
unchanged. It is easy to see that with (4.18) we have
$\omega_k$ unchanged but 
$$
\alpha_{j+1} \leftrightarrow \beta_j .
\leqno{(4.19)}
$$
Since with (4.18) we also have 
$$
U_n (x) \leftrightarrow V_n (x) 
\leqno{(4.20)}
$$
where
$$
V_n (x) = \W \left( \frac{a\qomN}{de}; be^{-\xi}, \mu b e^\xi ,
\frac{\qomN}{d}, \frac{\qomN}{e}, \frac{aq^2}{de}, 
\frac{a^2q^{1-N+n}}{\mu b^2 de} , \qmn ; q \right) .
\leqno{(4.21)}
$$
It follows from (4.17)--(4.21) that we have a full
biorthogonality
$$
\sum_{k=0}^N U_n (x_k) V_m (x_k) \omega_k =0, \qquad n
\neq m .
\leqno{(4.22)}
$$
It remains to calculate the $n=m$ case to obtain
$$
W_n \coloneq \sum_{k=0}^n U_n (x_k) V_n (x_k) \omega_k .
\leqno{(4.23)}
$$
This, the most tedious part of the calculation, we
now outline below. 

From Ismail and Masson [17] we have a fundamental
expression for the Stieltjes transform
$$
\eqalign{
& \sum_{k=0}^N {{Q_m (x_k)P_n (x_k) R_k} \over {\left[
\prod_{j=1}^n (x_k -\alpha_{j+1})(x_k - \beta_{j+1} )
\right] (x - x_k) }} \cr
& \qquad = {{Q_m (x) Y_n^{ (\min)}(x)} \over {\lambda_1 Y_{-1}
^{(\min)} (x) \left[ \prod_{j=0}^n (x-\alpha_{j+1})(x -
\beta_{j+1}) \right] }} , \qquad 0 \leq m < n \cr
}
\leqno{(4.24)}
$$
where $Y_n^{ (\min)}$ is the minimal solution to the
recurrence (4.3$'$) and $Q_m (x)$ is any polynomial of
degree $m <n$. The calculation of $W_n$ in (4.23) will
require an evaluation of the right side of (4.24) at $x =
\alpha_{n+2}$ for $m=n$. 

We first note that 
$$
{{P_n(x_k)} \over {\prod_{j=1}^n (x_k - \beta_{j+1})}} =
{{D_n U_n (x_k)} \over {\prod_{k=-1}^n u_k}}
\leqno{(4.25)}
$$
where
$$
D_n = \left( - \frac{2a}{b} \right)^n q^{n(n+1)}
{{\left(\frac{s}{q}, \frac{s}{aq}, \frac{aq}{d},
\frac{aq}{e}, \qmN \right)_n } \over { \left(\frac{s}{q}
\right)_{2n} (aq)_n}}  .
\leqno{(4.26)}
$$
Secondly, we may write
$$
V_n (x_k) = {{ E_n Q_n (x_k)} \over {\prod_{j=1}^n (x_k
-\alpha_{j+2} ) }} + {\rm\ additional\ terms}
\leqno{(4.27)}
$$
where 
$$
\eqalign{
E_n & = q^{{-n^2}} \left( \frac{aq}{2\mu bs} \right)^n
{{\left( 1 - \frac{aq^{2n-N+1}}{de} \right)} \over { \left(
1 - \frac{a\qmNpo}{de} \right)}}
{{\left( 
\frac{\qmNpo}{d}, 
\frac{\qmNpo}{e}, 
\frac{aq^2}{de}, 
sq^{n-1}, 
\frac{a\qmNpo}{de} \right)_n} \over {\left(
\frac{aq}{d}, 
\frac{aq}{e}, 
\qmN, 
\frac{aq^{3-N-n}}{des},
\frac{aq^{2+n-N}}{de} \right)_n }} , \cr
Q_n (x) & = (be^{-\xi} , \mu b e^\xi)_n . \cr
}
\leqno{(4.28)}
$$
Note that the additional terms in (4.27) do not
contribute to the right side of (4.24) evaluated at $x=
\alpha_{n+2}$ because they have a numerator factor $(x -
\alpha_{n+2})$. From (4.24)--(4.28) it now follows that
(4.23) can be calculated since 
$$
W_n = - {{E_n F_n G_n (1-a\qmNpo)\left(1-\mu b^2
\frac{de}{a} \qNmo \right) u_0 } \over {2 D_n H_n R_0 \mu
b}}
\leqno{(4.29)}
$$
where $E_n$ and $D_n$ are given by (4.26) and (4.28),
$$
R_0 = - {{u_0 aq(1-aq)\left( aq^2 , \frac{\mu b^2 d}{a},
\frac{\mu b^2 e}{a}, \frac{de}{a} \right)_N} \over {2\mu
bs (1-d)(1-e)(1-a\qNpo)\left( \mu b^2 q ,dq, eq,
\frac{\mu b^2 de}{a^2 q} \right)_N }} ,
\leqno{(4.30)}
$$
$$
G_n = (be^{-\xi}, \mu be^\xi)_n \vert_{x = \alpha_{n+2}}
= \left( \frac{a}{s} \qomn , \frac{\mu b^2 s}{a} \qnmo
\right)_n , 
\leqno{(4.31)}
$$
$$
\eqalign{
H_n & = \prod_{k=1}^n (x - \alpha_{k+1} )(x - \beta_{k+1}
) \vert_{x = \alpha_{n+2}} \cr
& = \left( \frac{q^2}{4\mu s}
\right)^n q^{-n(n+1)}
\left( \qmn , \frac{\mu b^2 s^2}{a^2} \qn
, s\qn , \frac{a^2q^{-n+2}}{\mu b^2 s} \right)_n , \cr
}
\leqno{(4.32)}
$$
and
$$
F_n = \left. {{X_n^{(\min)} (x)} \over {\gamma_0 X_{-1}^{(\min)}
(x) (x - \alpha_1)(x-\beta_1)}} \right\vert_{x =
\alpha_{n+2}} .
\leqno{(4.33)}
$$
From the minimal solution expression (3.30), with the
replacement (4.1), we find that
$$
\eqalign{
F_n & = {{ q^{-1} (1-s\qtnmo)(1-aq)\left( \frac{q^2 a}{s}
\right)^{n+1} } \over { \left( 1 - \frac{s}{q} \right) 
\left( 1 - \frac{a}{s} \qomn \right)
\left( 1 - \mu\frac{b^2}{a} s \qnmo \right)
(1-d)(1-e) (1-a\qNpo)}} \cr
& \qquad \times {{\left( q, \qmn , \frac{\mu b^2 s^2}{a^2}
\qnmt , \frac{ds}{aq}, \frac{es}{aq} , s\qN, s\qn ,
\frac{a^2\qtmn}{\mu b^2 s} , \frac{aq}{d} , \frac{aq}{e},
\qmN , \frac{s}{q} \right)_n (aq^2)_{2n} } \over {
(s,s)_{2n} \left( \frac{a\qtmn}{s} , \frac{\mu b^2 s\qn}{a}
, dq, eq, a\qNpt , aq \right)_n }} \cr
& \qquad \times \w \left( a\qtnpo ; \qnpo,  \frac{a^2q^2}{\mu
b^2 s} , \frac{a\qnpo}{d} , \frac{a\qnpo}{e}, \qmNpn ; q
\right) . \cr
}
\leqno{(4.33)}
$$
Using the Jackson summation for this last $\w$ and
putting all factors together we finally obtain, after much
simplification,
$$
W_n = \qmn {{\left( q, 
\frac{a\qomN}{\mu b^2 d}, 
\frac{a\qomN}{\mu b^2 e} ,
\frac{a^2q^2}{\mu b^2 e} ,
\frac{a\qtmN}{de} , aq \right)_n } \over {\left( \qN, 
\frac{a\qomN}{\mu b^2 de}, 
\frac{aq}{d},
\frac{aq}{e},
\frac{a}{\mu b^2},
\frac{a^2\qtmN}{\mu b^2 de} \right)_n }}
{{\left( \mu b^2 q, d,e,\frac{\mu b^2
de}{a^2q} \right)_N \left( 1 - \frac{a^2q^{1-N+n}}{\mu b^2
de} \right) } \over { \left( aq, 
\frac{\mu b^2 d}{a}, 
\frac{\mu b^2 e}{a}, 
\frac{de}{aq} \right)_N \left( 1 - \frac{a^2 q^{1-N+2n}
}{\mu b^2 de} \right) }} .
\leqno{(4.34)}
$$
We now summarize all of the above but with a normalized
probability measure. That is, we choose
$$
r_k \coloneq {{\omega_k} \over {W_0}}
\leqno{(4.35)}
$$
so that 
$$
\sum_{k=0}^N r_k  = 1 . 
\leqno{(4.36)}
$$

\noindent{\bf Theorem 4.1.} {\it Let $x = {{\left( e^\xi +
\mu^{-1} e^{-\xi} \right)} \slash 2}$, \/}
$$
\eqalign{
U_n (x) & = \W \left( a ; be^{-\xi} , \mu be^\xi,d,e,a\qNpo
, \frac{a^2q^{1-N+n}}{\mu b^2 de} , \qmn ; q \right) ,
\cr
V_m (x) & = \W \left( \frac{a\qomN}{de} ; be^{-\xi} , \mu be^\xi,
\frac{\qomN}{d},
\frac{\qomN}{e},
\frac{aq^2}{de}
, \frac{a^2q^{1-N+m}}{\mu b^2 de} , q^{-m} ; q \right) .
\cr
}
$$
{\it Then \/}
$$
\sum_{k=0}^N U_n (x_k)V_m (x_k) r_k = C_n \delta_{n,m},
\qquad 0 \leq n, \ m \leq N
\leqno{(4.37)}
$$
{\it where \/}
$$
\eqalign{
x_k & = {{ \left(bq^k + {{q^{-k}} \slash {b\mu}} \right)}
\slash 2} , \cr
r_k & = q^k {{\left(
\mu b^2,
q\sqrt{\mu b^2},
- q\sqrt{\mu b^2},
\frac{\mu b^2}{a} ,
\frac{aq}{d} ,
\frac{aq}{e} ,
\frac{\mu b^2 de}{a} 
\qNmo , \qmN \right)_k } \over
{\left( 
\sqrt{\mu b^2},
- \sqrt{\mu b^2}, aq, 
\frac{\mu b^2 d}{a} , 
\frac{\mu b^2 e}{a} , 
\mu b^2 \qNpo, \frac{a\qtmN}{de} , q \right)_k }}
{{\left( aq, 
\frac{\mu b^2 d}{a},
\frac{\mu b^2 e}{a},
\frac{de}{aq} \right)_N} \over {\left( \mu
b^2q,d,e,\frac{\mu b^2 de}{a^2q} \right)_N}} , \cr
C_n & =\qmn {{\left( q,
\frac{a\qomN}{\mu b^2 d},
\frac{a\qomN}{\mu b^2 e},
\frac{a^2q^2}{\mu b^2 de},
\frac{a\qtmN}{de}, 
a\qN \right)_n} \over {\left( \qmN ,
\frac{a\qomN}{\mu b^2 de},
\frac{aq}{d},
\frac{aq}{e},
\frac{a}{\mu b^2},
\frac{a^2\qtmN}{\mu b^2 de} \right)_n }}
{{
\left( 1 - \frac{a^2q^{1-N+n}}{\mu b^2 de} \right)
} \over {
\left( 1 - \frac{a^2q^{1-N+2n}}{\mu b^2 de} \right) }} .
\cr
}
$$

We have previously stated Theorem 4.1 for the  special
case $\mu = -1$ [25]. Wilson derived Theorem 4.1 for the
special case $\mu =1$ [30] but with a misprint. (In $C_n$
the factor $\qmn$ was omitted and the numerator factor
$(aq)_N$ in $r_k$ was incorrectly written as $(a^2q)_N$.)
Rahman and Suslov [26] obtained a general biorthogonality
which we have not checked but which should be equivalent
to Theorem 4.1 and its $q \rightarrow 1$ limit given
in Corollary 4.7 below.

The two main ingredients in the derivation of Theorem 4.1 were
the continued fraction (3.34), which came from a three-term
recurrence
for ${}_{10}\phi_9$'s and the Jackson $_8\phi_7$ summation
formula, which is also the "$q$-beta integral" in (4.37) (the
$n=m=0$
case). This gives a $q$-version of Askey's conjecture [1, p. 37] of
a connection
between Ramanujan's Entry 40, Dougall's $_7F_6$ summation formula
and a
${}_9F_8$ three-term recurrence.

We now state six limit cases of Theorem 4.1 as Corollaries.
The first five are at the ${}_4\phi_3$ (or very-well-
poised $\ephis$) level while the final one is a $q
\rightarrow 1$ limit at the ${}_9 F_8$ level. Further
lower level limits may be taken. This will result in an
Askey type scheme for discrete rational biorthogonality
analogous to that for hypergeometric orthogonal
polynomials [20]. This new scheme of orthogonality makes
contact with the polynomial scheme through Corollary 4.6
which is the case of $q$-Racah polynomials.
\medskip
{\bf 4 b). The ${}_4\phi_3$ level.}
We give five limiting cases of Theorem 4.1. The
calculations are straightforward 
and will not be detailed.

\medskip
\noindent{\bf Corollary 4.2.} {\it Let\/}
$$
\eqalign{
U_n (x) & = {{\left(aq, \frac{aq}{de} \right)_n} \over
{\left( \frac{aq}{d}, \frac{aq}{e} \right)_n }}
{}_4\phi_3 \left( {{\frac{2x\qmN}{b} ,d,e,\qmn} \atop
{\frac{2aqx}{b} , \qmN , \frac{de\qmn}{a} }} ; q \right)
\cr
V_m (x) & = {{\left(\frac{a\qmNpt}{de}, a\qN \right)_m} \over
{\left( \frac{aq}{e}, \frac{aq}{d} \right)_m }}
{}_4\phi_3 \left( {{\frac{2x\qmN}{b} ,\frac{\qmNpo}{d},
\frac{\qmNpo}{e},q^{-m}} \atop
{\frac{2ax\qmNpt}{bde} , \qmN , \frac{q^{-N-m+1}}{a} }} ; q 
\right) .
\cr
}
$$
{\it Then\/}
$$
\sum_{k=0}^N U_n (x_k) V_m (x_k) r_k =C_n \delta_{n,m} ,
$$
{\it where\/} $x_k = {1 \over 2} b\qk$,
$$
r_k = {{\left( \frac{aq}{d} , \frac{aq}{e} ,
\qmN\right)_k } \over { \left( aq, \frac{a\qtmN}{de} , q
\right)_k }} \qk {{\left( aq , \frac{de}{aq} \right)_N }
\over { (d,e)_N}} ,
$$
{\it and\/}
$$
C_n = \qmn {{\left( q, aq, \frac{a\qtmN}{de} \right)_n}
\over { \left( \qmN , \frac{aq}{d} , \frac{aq}{e}
\right)_n }} .
$$

\noindent{\bf Proof:} Take the limit $\mu \rightarrow
\infty$ in Theorem 4.1. The limiting $U_n$ and $V_m$ are
terminating $\ephis$'s. We then use Watson's transformation
formula [7, ({\caprom{3}}.~18), p.~242] 
to obtain the above
${}_4\phi_3$ expressions for $U_n(x)$ and $V_m (x)$. 

\medskip
\noindent{\bf Corollary 4.3.} {\it Let\/} $x= {1 \over 2}
( e^\xi +\mu^{-1} e^{-\xi})$,
$$
\eqalign{
U_n (x) & = {{\left(aq, \frac{aq}{\mu be} e^{-\xi} \right)_n} \over
{\left( \frac{aq}{\mu b} e^{-\xi}, \frac{aq}{e} \right)_n }}
{}_4\phi_3 \left( {{\frac{aq}{bd}e^\xi ,\mu be^\xi,e,\qmn} \atop
{\frac{aq}{b}e^\xi , \frac{aq}{d}, \frac{\mu be}{a}
e^\xi \qmn}} ; q \right)
\cr
V_m (x) & = {}_4\phi_3 \left( {{be^{-\xi}, \mu be^\xi,
\frac{aq^2}{de}, q^{-m} }\atop
{\frac{aq}{e},\frac{aq}{d},\frac{\mu b^2}{a}q^{1-m}}}; 
q\right).\cr
}
$$
{\it Then\/}
$$
\sum_{k=0}^\infty U_n (x_k) V_m (x_k) r_k = C_n
\delta_{n,m} ,
$$
{\it where\/}
$$
\eqalign{
x_k & = {1 \over 2} \left( b\qk + \qmk \slash b\mu \right), \cr
r_k & = {{\left( \mu b^2 , q\sqrt{\mu b^2},
- q\sqrt{\mu b^2}, \frac{\mu b^2}{a}, \frac{aq}{d},
\frac{aq}{e} \right)_k} \over {\left( \sqrt{\mu b^2},
- \sqrt{\mu b^2}, aq, \frac{\mu b^2 d}{a},
\frac{\mu b^2 e}{a}, q \right)_k }} \left( \frac{de}{aq}
\right)^k {{\left( aq, \frac{\mu b^2 d}{a},
\frac{\mu b^2 e}{a},
\frac{de}{aq} \right)_\infty } \over { \left( \mu b^2
q,d,e, \frac{\mu b^2 de}{a^2 q} \right)_\infty }} \cr
}
$$
{\it and\/}
$$
C_n = {{ \left( q, \frac{a^2q^2}{\mu b^2 de},aq \right)_n
} \over { \left( \frac{aq}{d} , \frac{aq}{e} ,
\frac{a}{\mu b^2} \right)_n }} .
$$

\noindent{\bf Proof:} Take the limit $N \rightarrow
\infty$ in Theorem 4.1 and use Watson's transformation
formula. 

\medskip
\noindent{\bf Corollary 4.4.} {\it Let\/} $x = {1 \over 2}
(e^\xi + \mu^{-1} e^{-\xi} )$,
$$
\eqalign{
U_n (x) & = 
{1 \over {d^n}} 
{{\left(aq, \frac{aq}{\mu b^2 e} \right)_n} \over
{\left( \frac{aq}{d} , \frac{aq}{\mu b^2 de} \right)_n }}
{}_4\phi_3 \left( {{\frac{aq}{\mu b^2}, d, \frac{a^2
\qopn}{\mu b^2 de}, \qmn} \atop
{\frac{aq}{b}e^\xi , \frac{aq}{b\mu}e^{-\xi} , \frac{aq}{\mu
b^2 e} }} ; q \right) , \cr
V_m (x) & = 
\left( {e \over q} \right)^m 
{{\left(\frac{aq^2}{de}, \frac{aq}{\mu b^2 e} \right)_m} \over
{\left( \frac{aq}{d} , \frac{a}{\mu b^2} \right)_m }}
{}_4\phi_3 \left( { {
\frac{aq^2}{\mu b^2 de}, \frac{q}{e}, \frac{a^2q^{m+1}}{\mu
b^2de}, q^{-m} 
}\atop {
\frac{aq^2}{bde} e^\xi, \frac{aq^2}{\mu bde}e^{-\xi} ,
\frac{aq}{\mu
b^2 e}
}} ; q \right) . \cr
}
$$
{\it Then\/}
$$
\sum_{k=0}^\infty U_n (x_k) V_m (x_k) r_k = C_n \delta_{m,n}
,
$$
{\it where\/}
$$
\eqalign{
x_k & = {1 \over 2} \left( b\qk + {{q^{-k}} \over {b\mu}}
\right) , \cr
r_k & = {{ \left( \mu b^2 ,
q \sqrt{\mu b^2},
- q \sqrt{\mu b^2},
\frac{\mu b^2}{a} ,
\frac{aq}{d} ,
\frac{\mu b^2 de}{aq} 
\right)_k } \over { \left(
\sqrt{\mu b^2},
- \sqrt{\mu b^2},
aq,
\frac{\mu b^2 d}{a} ,
\frac{aq^2}{de} , q 
\right)_k }} \left( \frac{aq}{\mu b^2 e} \right)^k {{\left(
aq, 
\frac{\mu b^2 d}{a} , \frac{aq}{\mu b^2 e} , \frac{aq}{d}
\right)} \over {\left( \mu b^2 q, d, \frac{q}{e},
\frac{a^2q^2}{\mu b^2 de} \right) }} \cr
}
$$
{\it and\/}
$$
C_n = {{\left(q, 
\frac{aq}{\mu b^2 e} , \frac{aq^2}{de} , aq \right)_n }
\over { \left( \frac{aq}{\mu b^2 de} , \frac{aq}{d} ,
\frac{a}{\mu b^2} , \frac{a^2 q^2}{\mu b^2 de} \right)_n }}
\left( \frac{a}{\mu b^2 d} \right)^n {{\left( 1- \frac{a^2
\qopn}{\mu b^2 de} \right) } \over { \left( 1 - \frac{a^2
q^{1 +2n} }{\mu b^2 de} \right) }} .
$$

\noindent{\bf Proof:} Replace $e$ by $e\qmN$ in Theorem 4.1,
take the limit as $N \rightarrow \infty$ and use Watson's
transformation formula. 

\medskip
\noindent{\bf Corollary 4.5.} {\it Let\/}
$$
\eqalign{
U_n (x) & = {{\left( aq, \frac{2aqx}{bd} \right)_n} \over {
\left( \frac{aq}{d} , \frac{2aqx}{b} \right)_n }} {}_4\phi_3
\left( {{
\frac{b^2 \mu de}{a^2q} \qmn ,d,\frac{b}{2x}, \qmn } \atop {
\frac{b^2 \mu de}{a} \qNmn , \qmN , \frac{bd}{2ax}\qmn }} ;
q \right) \cr
V_m (x) & = {}_4\phi_3 \left( {{
\frac{b}{2x}, \frac{\qmNpo}{d}, \frac{a^2q^{1-N+m}}{b^2 \mu
de} , \qmm} \atop { \frac{a\qmNpt}{2\mu bdex} , \frac{aq}{d}, 
\qmN}} ; q \right). \cr
}
$$
{\it Then\/}
$$
\sum_{k=0}^N U_n (x_k) V_m (x_k) r_k = C_n \delta_{n,m} ,
$$
{\it where\/} $x_k = {1 \over 2} b\qk$,
$$
r_k = {{\left( \frac{aq}{d} , \qmN , \frac{\mu b^2 de}{aq}
\qN \right)_k } \over {\left( aq, \frac{\mu b^2e}{a} ,q
\right)_k }}\qk \left( \frac{d}{aq} \right)^N {{\left( aq ,
\frac{\mu b^2 e}{a} \right)_N} \over {\left( d,\frac{\mu b^2
de}{a^2 q} \right)_N }} ,
$$
{\it and\/}
$$
C_n = \left( \frac{\qmN}{d} \right)^n {{\left( q ,
\frac{a\qomN}{\mu b^2 e} , \frac{a^2q^2}{\mu b^2 de} , aq
\right)_n } \over { \left( \qmN , \frac{a\qmN}{\mu b^2 de} ,
\frac{aq}{d} , \frac{a^2 \qtmN}{\mu b^2 de} \right)_n }}
{{\left( 1 - \frac{a^2 q^{1-N+n}}{\mu b^2 de} \right)}
\over {\left( 1 - \frac{a^2q^{1-N+2n}}{\mu b^2 de} \right)}}
.
$$

\medskip
\noindent{\bf Proof:} Replace $e$ and $\mu$ by $e\qM$ and
$\mu \qmM$ respectively in Theorem 4.1, take the limit as $M
\rightarrow \infty$ and use Watson's transformation formula. 

\medskip
\noindent{\bf Corollary 4.6.} {\it Let\/} $x = {1 \over 2}
(e^\xi +\mu^{-1} e^{-\xi} )$,
$$
\eqalign{
U_n (x) & = {}_4\phi_3
\left( {{
be^{-\xi}, \mu be^\xi, \frac{a^2 q^{1-N+n}}{\mu b^2 de},
\qmn} \atop {
\frac{aq}{d} , \frac{aq}{e}, \qmN}} ;
q \right) \cr
V_m (x) & = {}_4\phi_3 \left( {{
be^{-\xi}, \mu be^\xi, \frac{a^2q^{1-N+m}}{\mu b^2 
de} , \qmm} \atop { \frac{aq}{d} , \frac{aq}{e}, 
\qmN}} ; q \right). \cr
}
$$
{\it Then\/}
$$
\sum_{k=0}^N U_n (x_k) V_m (x_k) r_k = C_n \delta_{n,m} ,
$$
{\it where\/} $x_k = {1 \over 2} \left( b\qk +
\frac{\qmk}{b\mu} \right)$,
$$
r_k = {{\left( 
\mu b^2 , q \sqrt{\mu b^2}, -q \sqrt{\mu b^2},
\frac{aq}{d} , \frac{aq}{e}, \qmN
\right)_k } \over {\left( 
\sqrt{\mu b^2} , - \sqrt{\mu b^2},
\frac{\mu b^2 d}{a} ,
\frac{\mu b^2 e}{a} , \mu b^2 \qNpo , q
\right)_k }}\qk \left( 
\frac{\mu b^2 de}{a^2} \qNmt \right)^k {{\left( 
\frac{\mu b^2 d}{a} , \frac{\mu b^2 e}{a} 
\right)_N} \over {\left( \mu b^2 q, \frac{\mu b^2de}{a^2q}
\right)_N }} ,
$$
{\it and\/}
$$
C_n = \left( \mu b^2 \right)^n {{\left( q ,
\frac{a\qomN}{\mu b^2 d} , \frac{a\qomN}{\mu b^2 e} , 
\frac{a^2q^2}{\mu b^2 de} 
\right)_n } \over { \left( \qmN , \frac{aq}{d} ,
\frac{aq}{e} , \frac{a^2 \qtmN}{\mu b^2 de} \right)_n }}
{{\left( 1 - \frac{a^2 q^{1-N+n}}{\mu b^2 de} \right)}
\over {\left( 1 - \frac{a^2q^{1-N+2n}}{\mu b^2 de} \right)}}
.
$$

\medskip
\noindent{\bf Proof:} Replace $a$, $d$, $e$ by $a\qmM$,
$d\qmM$, $e\qmM$ respectively in Theorem 4.1 and take 
the limit as $M \rightarrow \infty$. 

Note that Corollary 4.6 corresponds to the case of $q$-Racah
polynomials
where $U_n (x) = V_n (x)$ are polynomials in $x$ of degree $n$.
\medskip
{\bf 4 c). The top ${}_9F_8$ level.} The $q \rightarrow 1$ limit
of Theorem 4.1 yields
\medskip
\noindent{\bf Corollary 4.7} {\it Let\/} $x = u(u+\mu)$,

$$
\eqalign{
U_n (x) & = {}_9F_8 \left[ {{a, 1+{1 \over 2}a, b-u,
b+u+\mu , d,e,a+N+1, } \atop { {1
\over 2}a, a - b+1+u, a-b-\mu +1-u, a-d+1, a-e+1, -N,
}} \right. \cr
& \qquad\qquad\qquad \left. {{
2a+1-N-2b-d-e-\mu +n, -n } \atop {
2b+\nu +d+e-a+N-n, a+n+1}} ; 1 
{\vphantom{ {{1 \over 2} \atop {1 \over 2}}}}
\right] \cr
V_m (x) & = {}_9F_8 \left[ {{ a-d-e-N+1, 1 + {1 \over
2}(a-d-e-N+1), b-u, b+\mu +u,
} \atop { {1 \over 2}(a-d-e-N+1), a - b-d-e-N+2+u, 
}} \right.  \cr
& \qquad\qquad\qquad {{ 
 -N-d+1, -N-e+1, a-d-e+2,
 } \atop {
a-\mu-b -d -e -N +2 -u , a-e+1, a-d+1, -N, }} \cr
& \qquad\qquad\qquad \left. {{  2a+1-N-2b-d-e-\mu +m, -m } \atop {
2b-a+\mu -m+1, a-d-e-N +2+m }} ; 1  
{\vphantom{ {{1 \over 2} \atop {1 \over 2}}}}
\right] . \cr
}
$$
{\it Then\/}
$$
\sum_{k=0}^N U_n (x_k) V_m (x_k) r_k = C_n \delta_{n,m} ,
$$
{\it where\/}
$$
\eqalign{
& x_k = (b+k)(b+k+\mu ), \cr
& r_k = 
{{
{\displaystyle{ 
{{ \left( \mu +2b,1 + {1 \over 2} (\mu+2b), \mu +2b
-a, a-d+1,\right. \qquad} \atop {\left. \qquad a-e+1, -N, 
\mu +2b +d+e-a+N -1\right)_k}}
}}
} \over {
{\displaystyle{
{{\left( {1 \over 2}(\mu+2b), a+1, \mu +2b +d-a, 
\mu +2b +e-a, \right. \qquad} \atop {\left. \qquad\qquad \mu +2b 
+N+1, a-d-e+2-N, 1 \right)_k }}
}}
}} \cr
& \qquad\times {{\left( a+1, \mu+2b+d-a, \mu+2b+e-a,
d+e-a-1 \right)_N} \over {\left(\mu +2b+1, d,e,\mu +2b
+d+e-2a-1 \right)_N }} , \cr
}
$$
{\it and\/}
$$
\eqalign{
C_n & = {{\left( 1,a+1-N-\mu -2b -d, a+1-N-\mu
-2b-e\right)_n} \over {\left( -N,a+1-N-\mu-2b-d-e,a-d+1
\right)_n}} \cr
& \qquad \times {{\left( 2a+2-\mu-2b-d-e, a+2-N-d-e, a+1
\right)_n} \over { \left( a-e+1, a-\mu -2b, 2a+2-N-\mu
-2b-d-e\right)_n}} . \cr
}
$$

\noindent{\bf Proof:} In the orthogonality of Theorem
4.1, we replace 
$a$, $b$, $d$, $e$, $\mu$ 
by 
$q^a$, $q^b$, $q^d$, $q^e$, $q^\mu$  and take the limit
as $q \rightarrow 1$. Note that in the statement of
Corollary 4.7, $(a,b,\dots ,c)_k$ is now the usual
multiple shifted factorial and ${}_9F_8$ is a generalized
hypergeometric function (see Bailey [4]).

For the $\mu =0$ case of Corollary 4.7 see Wilson [30].

\bigskip
\noindent{\bf 5. An ${_8 \phi_7}$ model}
\medskip

If we replace $f$ by $Fq^M$ and $s$ by $\frac{a^3q^3q^{-M}}{bcdeF}$
in equation (3.2) and let $M \rightarrow \infty$, then (3.2)
reduces 
to the equation
$$
\eqalign{
& Y_{n+1} - c_n Y_n + d_n Y_{n-1} =0 , \cr
& c_n = C_n + D_n + \frac{a^2\qtnpt}{bcdeh^2} 
\frac{(1-b)(1-c)(1-d)(1-e)}
{\left(1-\frac{a\qnp}{h}\right)} , \cr
& d_n = C_{n-1}D_n \cr
& C_n = - {{
\left(1- \frac{a\qnp}{bh}\right)
\left(1- \frac{a\qnp}{ch}\right)
\left(1- \frac{a\qnp}{dh}\right)
\left(1- \frac{a\qnp}{eh}\right)
} \over {
\left(1- \frac{a\qnp}{h}\right) }} \cr
& D_n = - q 
\left(1- \frac{\qn}{h}\right)
\left(1- \frac{a\qn}{h}\right)
\left(1- \frac{a^2\qnp}{bcdeh}\right) . \cr
}
\leqno{(5.1)}
$$
In [12], (5.1) was derived from a three-term contiguous relation
for a
very-well-poised ${_8\phi_7}$ series in the particular case $h=1$.
Following the calculation in [12] we have the following four
solutions
of (5.1):
$$
\leqalignno{
Y_n^{(1)} & = {{(-1)^n \left( \frac{a\qnp}{h}\right)_\infty 
} \over { ( \frac{a\qnp}{bh}, \frac{a\qnp}{ch}, \frac{a\qnp}{dh},
\frac{a\qnp}{eh})_\infty }} {\w} \left(a;b,c,d,e,hq^{-n};
\frac{a^2\qnpt}{bcdeh} \right), &
{(5.2)} \cr
Y_n^{(2)} & = {{(-1)^n  } \over { ( \frac{\qnp}{h},
\frac{a^2\qnpt}{bcdeh}, \frac{a\qn}{h}
)_\infty }} {\w} \left(\frac{q}{a};\frac{q}{b},\frac{q}{c},
\frac{q}{d},\frac{q}{e}, \frac{\qnp}{h}; \frac{bcdeh}{a^2q} 
q^{-n} \right), &
{(5.3)} \cr
Y_n^{(3)} & = {{(-1)^n \left( \frac{a\qnp}{h}\right)_\infty 
} \over { ( \frac{a\qnp}{bh}, \frac{a\qnp}{ch}, \frac{a\qnp}{dh},
\frac{a\qnp}{eh})_\infty }} {\w} \left( \frac{bcde}{aq}; 
b,c,d,e, \frac{bcdehq^{-n-1}}{a^2}; \frac{\qnp}{h} \right), &
{(5.4)} \cr
{\rm and} \qquad  &  \cr
Y_n^{(4)} & = {{(-1)^n  } \over { ( \frac{\qnp}{h},
\frac{a^2\qnpt}{bcdeh}, \frac{a\qn}{h}
)_\infty }} {\w} \left(\frac{aq^2}{bcde};\frac{q}{b},\frac{q}{c},
\frac{q}{d},\frac{q}{e}, \frac{a^2\qnpt}{bcdeh}; h 
q^{-n} \right). &
{(5.5)} \cr
}
$$
In addition to the above four solutions, two more solutions may be
obtained as follows. In (5.1) we make the parameter replacements
$$
(a,b,c,d,e,h) \rightarrow \left(
\frac{B^2}{A}, B,
\frac{BC}{A},
\frac{BD}{A},
\frac{BE}{A},
\frac{BH}{A} \right)
$$
and renormalize so as to arrive at (5.1) with lower case letters 
$a$, $b$, $c$, $d$, $e$, $h$ replaced by capitals. Consequently 
we obtain a fifth solution to (5.1) viz., 
$$
Y_n^{(5)} = {{(-1)^n \left( \frac{b\qnp}{h}\right)_\infty
} \over { ( 
\frac{\qnp}{h},
\frac{a\qnp}{ch},
\frac{a\qnp}{dh},
\frac{a\qnp}{eh})_\infty }} {\w} \left(
\frac{b^2}{a}; 
b,\frac{bc}{a},\frac{bd}{a},\frac{be}{a}, \frac{{bhq^{-n}}}{a}; 
\frac{a^2\qnpt}{bcdeh} 
\right) .
\leqno{(5.6)}
$$
The ``reflection transformation'' (see [12]) applied to (5.1) and
(5.6) yields 
another solution 
$$
Y_n^{(6)} = {{(-1)^n  
} \over { ( 
\frac{a\qnp}{bh},
\frac{b\qn}{h},
\frac{a^2\qnpt}{bcdeh}
)_\infty }} {\w} \left(\frac{aq}{b^2};\frac{q}{b},\frac{aq}{bc},
\frac{aq}{bd},\frac{aq}{be}, \frac{a\qnp}{bh};
\frac{bcdehq^{-n-1}}{a^2} 
\right).
\leqno{(5.7)}
$$
Parameter interchanges $b \leftrightarrow (c,d,e)$ in (5.6) and
(5.7) give six more solutions. Thus we have twelve pairwise
linearly
independent solutions to the three-term recurrence (5.1). 

All
these solutions may be derived as limiting cases of our $\Phi$
solutions obtained in Section 3. Thus it is easily seen that 
$Y_n^{(1)}$, $Y_n^{(2)}$, $Y_n^{(3)}$, $Y_n^{(4)}$, $Y_n^{(5)}$ and
$Y_n^{(6)}$ are limiting cases of the solutions
$X_n^{(1), \frac{s}{h} \qnm}$, $X_n^{(2), \frac{h}{s} \qmnpt}$,
$X_n^{(4), \frac{aq^2}{s}}$, $X_n^{(3), \frac{s}{aq}}$,
$X_n^{(3), \frac{bs}{ah} \qnm}$ and 
$X_n^{(4), \frac{ah}{bs} \qmnpt}$ respectively. Limits of all the 
remaining $X_n$ solutions are either one of the $Y_n$ solutions 
or a linear combination of the $Y_n$ solutions e.g., it can be
shown 
that 
$X_n^{(1), b}$ gives a linear combination of $(Y_n^{(1)},
Y_n^{(5)})$;
$X_n^{(1), hq^{-n}}$ of $(Y_n^{(1)}, Y_n^{(4)})$;
$X_n^{(2), \frac{q}{b}}$ of $(Y_n^{(2)}, Y_n^{(6)})$;
$X_n^{(2), \frac{\qnp}{h}}$ of $(Y_n^{(2)}, Y_n^{(3)})$;
$X_n^{(5), \frac{a\qnp}{bh}}$ of $(Y_n^{(3)}, Y_n^{(6)})$;
$X_n^{(6), \frac{bh}{a}q^{-n}}$ of $(Y_n^{(4)}, Y_n^{(5)})$;
$X_n^{(7), \frac{bc}{a}}$ of $(Y_n^{(5)}$, a $b \leftrightarrow c$
interchange of $Y_n^{(5)})$;
$X_n^{(8), \frac{aq}{bc}}$ of $(Y_n^{(6)}$, a $b \leftrightarrow c$
interchange of $Y_n^{(6)})$ etc.

Any three of the twelve solutions of (5.1) are connected
by a three-term ${_8 \phi_7}$ transformation formula. It
can be shown that the connection is provided either
by the standard three-term ${_8 \phi_7}$ formula ([7],
({\caprom{3}} .~37), p.~246) or its iterate which we state
below explicitly in a form suited to our purpose:
$$
\eqalign{
& {\w} \left(a;b,c,d,e,f; \frac{a^2q^2}{bcdef} \right) \cr
& \qquad  = {{( aq, b,
\frac{b}{a}, \frac{bq}{a}, \frac{cq}{a}, \frac{dq}{a},
\frac{eq}{a},
\frac{fq}{a}, \frac{aq}{cd}, \frac{aq}{ce}, \frac{aq}{cf},
\frac{aq}{de}, \frac{aq}{df}, \frac{aq}{ef},
\frac{bcdef}{a^2q} )_\infty } \over { (
\frac{q^2}{a}, \frac{q}{c}, \frac{q}{d}, \frac{q}{e}, \frac{q}{f},
\frac{aq}{c}, \frac{aq}{d}, \frac{aq}{e}, \frac{aq}{f},
\frac{bc}{a},
\frac{bd}{a}, \frac{be}{a}, \frac{bf}{a}, \frac{cdef}{a^2},
\frac{a^2q}{cdef} )_\infty}} \cr
& \qquad  \qquad  \times {\w} \left(\frac{q}{a};\frac{q}{b},
\frac{q}{c},\frac{q}{d},
\frac{q}{e},\frac{q}{f}; \frac{bcdef}{a^2q}  \right)
+ {{( aq, c,
\frac{c}{a}, \frac{bq}{a}, \frac{bq}{d}, \frac{bq}{e},
\frac{bq}{f},
\frac{aq}{bd}, \frac{aq}{be}, \frac{aq}{bf}
)_\infty } \over { (
\frac{b^2q}{a}, \frac{q}{d}, \frac{q}{e}, \frac{q}{f},
\frac{c}{b}, \frac{aq}{b}, \frac{aq}{d}, \frac{aq}{e},
\frac{aq}{f}, \frac{bc}{a}
)_\infty}} \cr
& \qquad  \qquad  \times \left[ 1 - {{ (b,
\frac{q}{b}, \frac{a}{b}, \frac{bq}{a}, \frac{cd}{a},
\frac{aq}{cd},
\frac{ce}{a}, \frac{aq}{ce}, \frac{cfq}{a}, \frac{a}{cf},
\frac{bdef}{a^2},
\frac{a^2q}{bdef} )_\infty } \over {(c,
\frac{q}{c}, \frac{a}{c}, \frac{cq}{a}, \frac{bd}{a},
\frac{aq}{bd},
\frac{be}{a}, \frac{aq}{be}, \frac{bfq}{a}, \frac{a}{bf},
\frac{cdef}{a^2},
\frac{a^2q}{cdef} )_\infty } } \right] \cr
& \qquad  \qquad  \times {\w} \left(\frac{b^2}{a}; b, 
\frac{bc}{a},\frac{bd}{a},
\frac{be}{a}, \frac{bf}{a}; \frac{a^2q^2}{bcdef}  \right) . \cr
}
\leqno{(5.8)}
$$
We may derive
(5.8) from the standard ${_8 \phi_7}$ transformation
as follows. Starting from [7, (III.~37), p.~246] we 
first apply [7, (III.~23), p.~243] to one of the ${_8 \phi_7}$'s 
in that formula, replacing 
$$
\w \left( \frac{ef}{c}; \frac{aq}{bc}, \frac{aq}{cd},
\frac{ef}{a}, e, f; \frac{bd}{a} \right)
$$
by
$$
{{( \frac{efq}{c}, \frac{bcdef}{a^2q}, \frac{bq}{c},
\frac{dq}{c} )_\infty } \over { ( \frac{bef}{a},
\frac{def}{a}, \frac{aq^2}{c^2},
\frac{bd}{a})_\infty }} \w \left(
\frac{aq}{c^2}; \frac{aq}{cf}, \frac{aq}{ce}, \frac{q}{c},
\frac{aq}{cb}, \frac{aq}{cd}; \frac{bcdef}{a^2q} \right) .
$$
If we now iterate the formula in the modified form and make use of
Slater's infinite product identity [27], [7, p.~138],
we can reduce the result to (5.8). 

The standard three-term $\ephis$ formula easily provides a
connection
between the solutions $Y_n^{(1)}$, $Y_n^{(3)}$, $Y_n^{(5)}$ and we 
actually have 
$$
\eqalign{
& \left( 
\frac{b^2q}{a}, aq, 
\frac{a}{b},
\frac{bq}{a},
\frac{aq}{cd},
\frac{aq}{ce},
\frac{aq}{de},
\frac{bcd}{a},
\frac{bce}{a},
\frac{bde}{a},
\frac{cde}{a} \right)_\infty Y_n^{(3)} \cr
& \qquad  = \left( 
aq,  c,d,e,
\frac{bq}{a},
\frac{bq}{c},
\frac{bq}{d},
\frac{bq}{e},
\frac{bcde}{a},
\frac{bcde}{a^2},
\frac{a^2q}{bcde} \right)_\infty Y_n^{(5)} \cr
& \qquad  \qquad  - \frac{a}{b} \left( 
\frac{b^2q}{a},
\frac{aq}{b},
\frac{aq}{c},
\frac{aq}{d},
\frac{aq}{e},
\frac{bc}{a},
\frac{bd}{a},
\frac{be}{a},
\frac{cde}{a},
\frac{aq}{cde},
\frac{bcde}{a} \right)_\infty Y_n^{(1)} . \cr
}
\leqno{(5.9)}
$$
On the other hand, (5.8) is better suited to provide a connection
between the solutions $Y_n^{(1)}$, $Y_n^{(2)}$ and $Y_n^{(5)}$ 
and they are related by
$$
Y_n^{(1)}  = R Y_n^{(2)} + S Y_n^{(5)},
\leqno{(5.10)}
$$
where
$$
\eqalign{
R & = - \frac{bdeh}{a^2} {{( aq, b, 
\frac{c}{a}, \frac{bq}{a}, \frac{dq}{a}, \frac{eq}{a},
\frac{aq}{cd}, \frac{aq}{ce}, \frac{aq}{de}, \frac{a}{b},
\frac{bq}{a} )_\infty} \over {(
\frac{q}{c}, \frac{q}{d}, \frac{q}{e}, \frac{a}{c},
\frac{aq}{b}, \frac{aq}{d}, \frac{aq}{e}, \frac{bc}{a},
\frac{bd}{a}, \frac{be}{a},
\frac{q^2}{a} )_\infty }}
{{( 
\frac{bcdeh}{a^2}, \frac{a^2q}{bcdeh}, \frac{a}{qh},
\frac{hq^2}{a} )_\infty } \over { (
\frac{cdeh}{a^2}, \frac{a^2q}{cdeh}, \frac{a}{bh},
\frac{bhq}{a} )_\infty }} \cr 
S & =  {{( aq,c,
\frac{c}{a}, \frac{bq}{a}, \frac{bq}{d}, \frac{bq}{e},
\frac{aq}{bd},
\frac{aq}{be} )_\infty } \over { (
\frac{b^2q}{a}, \frac{q}{d}, \frac{q}{e}, \frac{c}{b},
\frac{aq}{b}, \frac{aq}{d}, \frac{aq}{e},
\frac{bc}{a} )_\infty }}
\left[ 1 - {{(b,
\frac{q}{b}, \frac{a}{b}, \frac{bq}{a}, \frac{cd}{a},
\frac{aq}{cd}, \frac{ce}{a}, \frac{aq}{ce}, \frac{chq}{a},
\frac{a}{ch}, \frac{bdeh}{a^2},
\frac{a^2q}{bdeh} )_\infty} \over {  (c,
\frac{q}{c}, \frac{a}{c}, \frac{cq}{a}, \frac{bd}{a},
\frac{aq}{bd}, \frac{be}{a}, \frac{aq}{be}, \frac{bhq}{a},
\frac{a}{bh}, \frac{cdeh}{a^2},
\frac{a^2q}{cdeh} )_\infty} } \right] . \cr
}
$$

We now give the continued fraction associated with the
three-term recurrence (5.1). We follow the same procedure
as in [12] where it was derived for the special case
$h=1$. For the sake of completeness, the outline of the
method and the related results obtained are being stated
below.

In order to construct the minimal solution to (5.1)
we examine the large $n$ asymptotics of the solutions
$Y_n^{(1)}$ and $Y_n^{(3)}$. For $Y_n^{(1)}$ we first
apply [7, ({\caprom {3}}.~23), p.~243] and then let $n
\rightarrow \infty$ to obtain
$$
\eqalign{
Y_n^{(1)} & \approx (-1)^n C_1, \qquad \left\vert
\frac{aq}{de} \right\vert < 1, \cr
C_1 & = {{\left( aq, \frac{aq}{de} \right)_\infty} \over
{ \left( \frac{aq}{d} , \frac{aq}{e} \right)_\infty }}
{}_3 \phi_2 \left( {{\frac{aq}{bc} ,d,e} \atop
{\frac{aq}{b} , \frac{aq}{c} }} ; \frac{aq}{de} \right) .
\cr
}
\leqno{(5.11)}
$$
For $Y_n^{(3)}$, we apply [7, ($\caprom 3$.~24), p.~243],
take the limit as $n \rightarrow \infty$ and subsequently
using [7, ($\caprom 3$.~9), p.~241] we have 
$$
\eqalign{
Y_n^{(3)} & \approx (-1)^n C_3, \qquad \left\vert
\frac{bc}{a} \right\vert < 1, \cr
C_3 & = {{\left( \frac{bcde}{a}, \frac{bc}{a}\right)_\infty} \over
{\left(\frac{bcd}{a} , \frac{bce}{a} \right)_\infty }}
{}_3 \phi_2 \left( {{\frac{de}{a} ,d,e} \atop
{\frac{cde}{a} , \frac{bde}{a} }} ; \frac{be}{a} \right) .
\cr
}
\leqno{(5.12)}
$$
A minimal solution to (5.1) is now given by 
$$
Y_n^{(\min)} = C_3Y_n^{(1)} - C_1 Y_n^{(3)} .
\leqno{(5.13)}
$$
Using Pincherle's theorem [8], [19] and simplifying we have the
continued fraction representation:

\medskip
\noindent{\bf Theorem 5.1.}
$$
\eqalign{
& \frac{1}{c_0} \lowminus \frac{d_1}{c_1} \lowminus
\frac{d_2}{c_2} \lowminus \lowdots
 = {{Y_0^{(\min)}} \over {c_0 Y_0^{(\min)} - Y_1^{(\min)}
}} \cr
& = - {1 \over {q \left(1 - \frac{1}{h}\right)\left( 1 -
\frac{a^2q}{bcdeh} \right) \left( 1 - \frac{a}{h} \right)
}} \cr 
& \times \left[ \vphantom{{{\frac{q}{q}} \over
{\frac{q}{q}}}} C_3 \, \w \left(a;b,c,d,e,h;
\frac{a^2q^2}{bcdeh} \right) \right. \cr
& \qquad  \left.  -C_1 {{\left( \frac{bcde}{a}, \frac{bcdeh}{a^2q},
\frac{bq}{h}, \frac{cq}{h}, \frac{dq}{h}, \frac{eq}{h}
\right)_\infty} \over { \left( \frac{bde}{a},
\frac{bce}{a}, \frac{bcd}{a}, \frac{cde}{a},
\frac{aq^2}{h^2}, \frac{q}{h} \right)_\infty }} \w
\left( \frac{aq}{h^2}; \frac{aq}{bh}, \frac{aq}{ch},
\frac{aq}{dh}, \frac{aq}{eh},
\frac{q}{h}; \frac{bcdeh}{a^2q} \right) \right] 
\cr
& \slash \left[  \vphantom{{{\frac{q}{q}} \over
{\frac{q}{q}}}} C_3 \, \w \left( a;b,c,d,e,hq;
\frac{a^2q}{bcdeh} \right) \right. \cr
& \qquad  \left.  - C_1 {{ \left( \frac{bcde}{a},
\frac{bcdeh}{a^2}, \frac{b}{h}, \frac{c}{h}, \frac{d}{h},
\frac{e}{h} \right)_\infty } \over {\left( 
\frac{bde}{a}, \frac{bce}{a}, \frac{bcd}{a}, \frac{cde}{a},
\frac{a}{h^2}, \frac{1}{h} \right)_\infty }} \w \left(
\frac{a}{h^2q}; \frac{a}{bh}, \frac{a}{ch}, \frac{a}{dh},
\frac{a}{eh}, \frac{1}{h}; \frac{bcdeh}{a^2} \right) \right]
.  \cr
}
\leqno{(5.14)}
$$
In the special case $h=1$, the above theorem gives (see [12]):

\medskip
\noindent{\bf Corollary 5.2.} {\it If $c_n$, $d_n$ are
given by (5.1) and $h=1$, we have the continued fraction
representation\/}
$$
\eqalign{
& \frac{1}{c_0} \lowminus \frac{d_1}{c_1} \lowminus
\frac{d_2}{c_2} \lowminus \lowdots \cr
& = \frac{bcde}{a^2q^2} \left\{ 
\vphantom{{{{{\frac{a}{a}}\atop{\frac{a}{a}}}}\over{{{\frac{q}{q}} 
\atop{\frac{q}{q}}}}}} 
{{(1-aq)} \over
{(1-b)(1-c)(1-d)(1-e)}} \w \left( aq; q, 
\frac{aq}{b},
\frac{aq}{c},
\frac{aq}{d},
\frac{aq}{e};
\frac{bcde}{a^2q} \right) \right. \cr
& \left.  - {{\left( \frac{bc}{a}, \frac{aq}{d}, \frac{aq}{e},
\frac{bde}{a}, \frac{cde}{a}, q \right)_\infty } \over {
\left( \frac{aq}{de}, b,c,d,e, \frac{bcde}{a^2q}
\right)_\infty }} {{  {}_3\phi_2 \left(
{{\frac{de}{a},d,e} \atop {\frac{cde}{a}, \frac{bde}{a}
}} ; \frac{bc}{a} \right) } \over { {}_3\phi_2 \left(
{{\frac{aq}{bc} ,d,e }\atop {\frac{aq}{b}, \frac{aq}{c}
}} ; \frac{aq}{de} \right) }} 
\vphantom{{{{{\frac{a}{a}}\atop{\frac{a}{a}}}}\over{{{\frac{q}{q}} 
\atop{\frac{q}{q}}}}}} 
\right\} . \cr
}
\leqno{(5.15)}
$$
On the other hand, if we write $f=\qmm$,
$s=\frac{a^3q^{3-m}}{bcde}$, $h=1$ in (3.32) and take the
limit as $m \rightarrow \infty$, we obtain the equivalent
continued fraction representation [10]:

\medskip
\noindent{\bf Corollary 5.2$'$.} {\it Under the same
conditions as Corollary 5.2 we have the equivalent
representation\/}
$$
\eqalign{
& \frac{1}{c_0} \lowminus \frac{d_1}{c_1} \lowminus
\frac{d_2}{c_2} \lowminus \lowdots \cr
& = {{\left( 1 - \frac{a}{q} \right) } \over { q
\left(1 - \frac{a}{b}\right)
\left(1 - \frac{a}{c}\right)
\left(1 - \frac{a}{d}\right)
\left(1 - \frac{a}{e}\right) }} \left[ \w \left( 
\frac{q}{a}; \frac{q}{b}, \frac{q}{c}, \frac{q}{d},
\frac{q}{e}, q; \frac{bcde}{a^2q} \right) - R \right] ,
\cr
R & = {{\left( q,a,\frac{q^2}{a}, \frac{de}{a},
\frac{dc}{a}, \frac{ec}{a} \right)_\infty} \over { \left(
\frac{dq}{a}, \frac{eq}{a}, \frac{cq}{a}, \frac{dec}{aq},
\frac{aq}{b}, b \right)_\infty }} \cr
& \times \left\{ {}_3\phi_2 \left( {{\frac{q}{d},
\frac{q}{e}, \frac{q}{c} } \atop {\frac{qb}{a},
\frac{aq^2}{cde} }} ; b \right) +  {{\left( 
\frac{q}{d}, \frac{q}{e}, \frac{q}{c}, \frac{bcde}{a^2},
\frac{aq}{b}, \frac{b}{a}, \frac{cde}{a^2}, \frac{a^2q}{cde},
\frac{cde}{aq} \right)_\infty } \over { \left(
\frac{ec}{a}, \frac{bq}{a}, \frac{aq}{cde}, \frac{q}{a}, a,
\frac{bcde}{a^2q}, \frac{a^2q^2}{bcde}, \frac{de}{a},
\frac{dc}{a} \right)_\infty}}
{{ {}_3\phi_2 \left( {{\frac{de}{a},
\frac{dc}{a}, \frac{ec}{a} } \atop {\frac{bcde}{a^2} ,
\frac{dec}{a} }} ; b \right) } \over {\left(
{{\frac{aq}{bc}, \frac{aq}{bd}, \frac{aq}{be} } \atop
{\frac{aq}{b}, \frac{a^2q^2}{bcde} }} ; b \right) }}
\right\} . \cr
}
\leqno{(5.16)}
$$
The above representation was obtained earlier in [10]. It
is not obvious that this is the same as the
representation (5.15). However, we can show that right
sides of (5.15) and (5.16) are indeed equal. The proof involves
use of the three-term $\ephis$ transformation formula
(5.8), ${}_3\phi_2$ transformation 
[7, ($\caprom 3$.~9), p.~241], 
three-term ${}_3\phi_2$ transformation formula
[7, ($\caprom 3$.~33), p.~245) and Bailey's infinite
product identity [7, Ex.~5.21, p.~138].

When the continued fraction terminates, we have a further
simplification (see [12]):

\medskip
\noindent{\bf Corollary 5.3.} {\it Under the conditions of
Corollary (5.2) but with one of $ \frac{aq}{b}, \frac{aq}{c},
\frac{aq}{d}, \frac{aq}{e} = \qmN, N=0,1,\dots$ we have\/}
$$
\eqalign{
& \frac{1}{c_0} \lowminus \frac{d_1}{c_1} \lowminus
\frac{d_2}{c_2} \lowminus \lowdots \lowminus \frac{d_N}{c_N} \cr
& = \frac{bcde}{a^2q^2} {{ (1-aq)} \over
{(1-b)(1-c)(1-d)(1-e)}} \w \left(aq; q, 
\frac{aq}{b},
\frac{aq}{c},
\frac{aq}{d},
\frac{aq}{e};
\frac{bcde}{a^2q} \right) . \cr
}
\leqno{(5.17)}
$$

It is this last continued fraction result which is
associated with the rational biorthogonality given in
Corollary 4.3 for the special case $e = a\qNpo$ [12].

\bigskip
\noindent{\bf 6. Three-term $\Phi$ transformation formula}
\medskip

Any three solutions out of the fifty-six solutions obtained 
for the second-order finite difference equation (3.2) in
Section 3 are connected by a three-term $\Phi$ transformation
formula. In order to derive  such a formula, we start by picking 
any three of the solutions and assume 
a linear dependence. Here we consider
$$
X_n^{(1), \frac{s}{h}\qnm} = P
X_n^{(2), \frac{\qnp}{h}} + Q 
X_n^{(3), \frac{bs}{ah}\qnm}
\leqno{(6.1)}
$$
where $P$ and $Q$ are independent of $n$. If we take the large-$n$
asymptotics of (6.1) then from (3.20), (3.23) and (3.24) we
obtain
$$
\eqalign{
& \w \left(a;b,c,d,e,f; \frac{s}{aq} \right) \cr
&\qquad  = P \w \left(
\frac{q}{a};
\frac{q}{b},
\frac{q}{c},
\frac{q}{d},
\frac{q}{e},
\frac{q}{f};
\frac{aq^2}{s} \right) \cr
&\qquad \qquad  + Q {{(
\frac{s}{aq},
\frac{bq}{c},
\frac{bq}{d},
\frac{bq}{e},
\frac{bq}{f},
\frac{bq}{a},
\frac{bhq}{s},
\frac{s}{bh} )_\infty } \over { ( 
\frac{b^2q}{a},
\frac{bc}{a},
\frac{bd}{a},
\frac{be}{a},
\frac{bf}{a}, b,
\frac{bh}{a},
\frac{aq}{bh} )_\infty }} \w \left(
\frac{b^2}{a}; b,
\frac{bc}{a},
\frac{bd}{a},
\frac{be}{a},
\frac{bf}{a};
\frac{s}{aq}  \right) . \cr
} 
\leqno{(6.2)}
$$
Therefore from (5.8)
$$
P = {{( aq,b, \frac{b}{a}, \frac{bq}{a}, \frac{cq}{a},
\frac{dq}{a}, \frac{eq}{a}, \frac{fq}{a}, \frac{aq}{cd},
\frac{aq}{ce}, \frac{aq}{cf} )_\infty } \over {(
\frac{q^2}{a}, \frac{q}{c}, \frac{q}{d}, \frac{q}{e},
\frac{q}{f}, \frac{aq}{c}, \frac{aq}{d}, \frac{aq}{e},
\frac{aq}{f}, \frac{bc}{a}, \frac{bd}{a})_\infty }}
{{( \frac{aq}{de}, \frac{aq}{df}, \frac{aq}{ef},
\frac{bcdef}{a^2q} )_\infty } \over {( \frac{be}{a},
\frac{bf}{a}, \frac{cdef}{a^2}, \frac{a^2q}{cdef} )_\infty }}
\leqno{(6.3)}
$$
and 
$$
\eqalign{
Q & = {{( aq,b,c,
\frac{c}{a},
\frac{bd}{a},
\frac{aq}{bd},
\frac{be}{a},
\frac{aq}{be},
\frac{bf}{a},
\frac{aq}{bf},
\frac{bh}{a},
\frac{aq}{bh} )_\infty } \over {(
\frac{s}{aq},
\frac{q}{d},
\frac{q}{e},
\frac{q}{f},
\frac{aq}{b},
\frac{aq}{d},
\frac{aq}{e},
\frac{aq}{f},
\frac{c}{b},
\frac{bq}{c},
\frac{s}{bh},
\frac{bhq}{s} )_\infty }} \cr
&\qquad  \times \left[ 1 - {{(b, 
\frac{q}{b},
\frac{a}{b},
\frac{bq}{a},
\frac{cd}{a},
\frac{aq}{cd},
\frac{ce}{a},
\frac{aq}{ce},
\frac{cfq}{a},
\frac{a}{cf},
\frac{bdef}{a^2},
\frac{a^2q}{bdef} )_\infty } \over { (c,
\frac{q}{c},
\frac{a}{c},
\frac{cq}{a},
\frac{bd}{a},
\frac{aq}{bd},
\frac{be}{a},
\frac{aq}{be},
\frac{bfq}{a},
\frac{a}{bf},
\frac{cdef}{a^2},
\frac{a^2q}{cdef} )_\infty }} \right] . \cr
}
\leqno{(6.4)}
$$
We now substitute these values of $P$ and $Q$ in (6.1) and 
replace $\left( a,b,c,d,e,f, \frac{s}{h} \qnm, h\qmn \right)$
by $(A,B,C,D,E,F,G,H)$ respectively. The result can be 
written as a three-term $\Phi$ transformation formula  
$$
\eqalign{
& \Phi^{(G)} (A;B,C,D,E,F,G,H;q) \cr
&\qquad  = {{( Aq,B,G,
\frac{B}{A}, \frac{G}{A}, \frac{Bq}{A},
\frac{Cq}{A}, \frac{Dq}{A}, \frac{Eq}{A},
\frac{Fq}{A},
\frac{Gq}{A} )_\infty } \over { (
\frac{q^2}{A}, \frac{q}{C}, \frac{q}{D},
\frac{q}{E}, \frac{q}{F}, \frac{q}{H},
\frac{A}{H}, \frac{Aq}{C}, \frac{Aq}{D},
\frac{Aq}{E}, \frac{Aq}{F},
\frac{Aq}{H} )_\infty }} \cr
&\qquad\qquad \qquad  \times {{(
\frac{Aq}{CD}, \frac{Aq}{CE}, \frac{Aq}{CF},
\frac{Aq}{DE}, \frac{Aq}{DF}, \frac{Aq}{EF},
\frac{Aq}{BH}, \frac{Aq}{CH}, \frac{Aq}{DH},
\frac{Aq}{EH}, \frac{Aq}{FH},
\frac{Aq}{GH} )_\infty } \over { (
\frac{BC}{A}, \frac{BD}{A}, \frac{BE}{A},
\frac{BF}{A}, \frac{BG}{A}, \frac{CG}{A},
\frac{DG}{A}, \frac{EG}{A}, \frac{FG}{A},
\frac{CDEF}{A^2},
\frac{A^2q}{CDEF} )_\infty }} \cr
&\qquad\qquad \qquad  \times \Phi^{(\frac{q}{H})} \left(
\frac{q}{A}; \frac{q}{B}, \frac{q}{C},
\frac{q}{D}, \frac{q}{E}, \frac{q}{F},
\frac{q}{G},
\frac{q}{H}; q \right) \cr
&\qquad \qquad  + {{(Aq, B,C,
\frac{C}{A}, \frac{G}{A}, \frac{Gq}{A},
\frac{Gq}{C}, \frac{Gq}{D}, \frac{Gq}{E},
\frac{Gq}{F}, \frac{Gq}{H} )_\infty } \over {(
\frac{G^2q}{A}, \frac{q}{D}, \frac{q}{E},
\frac{q}{F}, \frac{q}{H}, \frac{Aq}{B},
\frac{Aq}{D}, \frac{Aq}{E}, \frac{Aq}{F},
\frac{Aq}{H},
\frac{B}{G} )_\infty }} \cr
&\qquad\qquad \qquad  \times {{(
\frac{BD}{A}, \frac{Aq}{BD}, \frac{BE}{A},
\frac{Aq}{BE}, \frac{BF}{A}, \frac{Aq}{BF},
\frac{BH}{A}, \frac{Aq}{BH} )_\infty } \over {(
\frac{C}{B}, \frac{Bq}{C}, \frac{GB}{A},
\frac{GC}{A}, \frac{GD}{A}, \frac{GE}{A},
\frac{GF}{A}, \frac{GH}{A} )_\infty }} \cr
&\qquad\qquad \qquad  \times \left[ 1- {{( B, 
\frac{q}{B}, \frac{A}{B}, \frac{Bq}{A},
\frac{CD}{A}, \frac{Aq}{CD}, \frac{CE}{A},
\frac{Aq}{CE}, \frac{CFq}{A}, \frac{A}{CF}, \frac{BDEF}{A^2},
\frac{A^2q}{BDEF} )_\infty } \over {(C, 
\frac{q}{C}, \frac{A}{C}, \frac{Cq}{A}, \frac{BD}{A},
\frac{Aq}{BD}, \frac{BE}{A}, \frac{Aq}{BE}, \frac{BFq}{A},
\frac{A}{BF}, \frac{CDEF}{A^2},
\frac{A^2q}{CDEF} )_\infty }} \right] \cr
&\qquad\qquad \qquad  \times \Phi^{(\frac{BG}{A})} \left(
\frac{G^2}{A}; \frac{GB}{A}, \frac{GC}{A}, \frac{GD}{A},
\frac{GE}{A}, \frac{GF}{A}, G,
\frac{GH}{A}; q \right) . \cr
}
\leqno{(6.5)}
$$

\bigskip
\noindent{\bf Particular cases }
\medskip

\item{1.} In the general formula, if we write $B=q^{-n}$, 
then it reduces to the terminating $\tphia$ transformation
\par
$$
\eqalign{
& \W (A;C,D,E,F,G,H,q^{-n};q) \cr
& = {{\left(Aq, C,G,
\frac{C}{A},
\frac{G}{A},
\frac{A\qnp}{D},
\frac{A\qnp}{E},
\frac{A\qnp}{F},
\frac{A\qnp}{H},
\frac{\qmnp}{D},
\frac{\qmnp}{E},
\frac{\qmnp}{F},
\frac{\qmnp}{H},
\frac{\qmnp}{A} \right)_\infty } \over {\left(
\frac{\qmtnp}{A},
\frac{q}{D},
\frac{q}{E},
\frac{q}{F},
\frac{q}{H},
\frac{Aq}{D},
\frac{Aq}{E},
\frac{Aq}{F},
\frac{Aq}{H}, A\qnp, C\qn, G\qn, 
\frac{C}{A}\qmn,
\frac{G}{A}\qmn \right)_\infty }} \cr
&\qquad  \times \W \left(
\frac{\qmtn}{A};
\frac{C\qmn}{A},
\frac{D\qmn}{A},
\frac{E\qmn}{A},
\frac{F\qmn}{A},
\frac{G\qmn}{A},
\frac{H\qmn}{A}, \qmn ; q  \right) . 
}
\leqno{(6.6)}
$$
The above formula can be derived also with the help of 
known formulas. Refer to the terminating $\tphia$ transformation 
formula [7, Ex. 2.19, p.~53] 
to which, if we apply the $\tphia$ transformation formula 
[7, Ex. 2.30, p.~56],  we obtain (6.6). 

\medskip
\item{2.} Writing $G= \qmn$ in (6.5) we obtain another relation
between two terminating $\tphia$'s as in the above paragraph. 
For the particular value $n=0$, the relation reduces to Slater's 
infinite product identity ([27], also see 
[7, Ex. 5.22, p.~138]).

\medskip
\item{3.} In the general formula (6.5) we interchange
the parameters $B \leftrightarrow C$ and then eliminate 
$\Phi^{(\frac{q}{H})} \left( \frac{q}{A}; \frac{q}{B}, \frac{q}{C},
\frac{q}{D}, \frac{q}{E}, \frac{q}{F}, \frac{q}{G},
\frac{q}{H}; q \right )$ from the two equations. Using Slater's 
infinite product identity we can write the result in the 
following form: 
\par

$$
\eqalign{
& {{(
\frac{Aq}{D},
\frac{Aq}{E},
\frac{Aq}{F},
\frac{Aq}{H},
\frac{q}{D},
\frac{q}{E},
\frac{q}{F},
\frac{q}{H} )_\infty } \over { \left( Aq,B,C,G, 
\frac{B}{A},
\frac{C}{A},
\frac{G}{A} \right)_\infty }} \W (A;B,C,\dots ,H;q) \cr
& = - \frac{A}{B} {{(
\frac{Bq}{D},
\frac{Bq}{E},
\frac{Bq}{F},
\frac{Bq}{H},
\frac{Aq}{BD},
\frac{Aq}{BE},
\frac{Aq}{BF},
\frac{Aq}{BH} )_\infty } \over {(
\frac{B^2q}{A}, B,
\frac{BC}{A},
\frac{BG}{A},
\frac{A}{B},
\frac{C}{B},
\frac{G}{B} )_\infty }} \W \left(
\frac{B^2}{A}; B, 
\frac{BC}{A},
\frac{BD}{A},
\dots 
\frac{BH}{A}; q \right) \cr
&\qquad  + {\rm \ idem \ } (B;C,G). \cr
}
\leqno{(6.7)}
$$
The above is a particular case of [7, Ex. 4.6, p.~122]
when the balance condition is satisfied. Slater has obtained
(6.7) from other considerations and then used it to derive her 
infinite product identity. However, we have shown here that, 
from our general formula, we can derive Slater's infinite 
product identity and consequently also (6.7). 

\medskip
\item{4.} In (6.5), if we take the value $\frac{BG}{A} =1$ and 
use Slater's infinite product identity we can reduce the result 
to Bailey's non-terminating extension of Jackson's $\ephis$ sum
[7, (II.~25), p.~238].
\par 

\bigskip
\noindent{\bf Limiting cases}
\medskip 

All the known $\ephis$ two-term and three-term transformation 
formulas may be obtained as limiting cases of our general
transformation formula (6.5). When we let one of the parameters
in the $\Phi (A;B,C,D,E,F,G,H;q)$ tend to $0$ and another 
tend to infinity so that the balance condition remains intact, 
we obtain in general a relation connecting three, four or five 
$\ephis$'s depending on the choice of parameters. The $\ephis$ 
series involved are 
$$
\w \left(A;B,C,D,E,F; \frac{A^2q^2}{BCDEF} \right),
$$
its reflection
$$
\w \left(\frac{q}{A};\frac{q}{B},\frac{q}{C},\frac{q}{D},
\frac{q}{E},\frac{q}{F}; \frac{BCDEF}{A^2q} \right)
$$
and complements of both, i.e., series of the type
$$
\w \left(\frac{B^2}{A};\frac{BC}{A},\frac{BD}{A},\frac{BE}{A},
\frac{BF}{A},B; \frac{A^2q^2}{BCDEF} \right)
$$
and
$$
\w \left(\frac{Aq}{B^2};\frac{Aq}{BC},\frac{Aq}{BD},\frac{Aq}{BE},
\frac{Aq}{BF},\frac{q}{B}; \frac{BCDEF}{A^2q} \right).
$$
We give below some of the interesting limiting cases.

\medskip
\item{1.} Replace $B$ and $G$ by $bq^{-m}$ and $gq^m$
respectively and let $m \rightarrow\infty$. We arrive
at a relation connecting five $\ephis$'s. However,
if we write $b=1$ in this relation after using the
balance condition to eliminate $g$, we just obtain the
transformation formula [7, ($\caprom 3$.~24), p.~243]
of a very-well-poised $\ephis$ series into another
very-well-poised $\ephis$-series.
\par 

\medskip
\item{2.} Replacing $B$ and $H$ by $bq^{-m}$ and $hq^m$ 
respectively and letting $m \rightarrow \infty$ leads in
general to a representation of an $\ephis$ as sum of a ${_4\phi_3}$
and a ${_4\psi_4}$. Using now the balance condition to eliminate
$h$ and then writing the special value $b=1$ yields the 
following representation of an $\ephis$ as a sum of 
two ${_4\phi_3}$'s:
$$
\eqalign{
& \w \left(A;C,D,E,F,G; \frac{A^2q^2}{CDEFG} \right) \cr
&\qquad  =  {{(Aq,E,
\frac{Aq}{CG},
\frac{Aq}{DG},
\frac{Aq}{FG},
\frac{A^2q^2}{CDEF} )_\infty } \over {(
\frac{Aq}{C},
\frac{Aq}{D},
\frac{Aq}{F},
\frac{Aq}{G},
\frac{E}{G},
\frac{A^2q^2}{CDEFG} )_\infty }} {_4\phi_3} \left( {{
\frac{Aq}{CE},
\frac{Aq}{DE},
\frac{Aq}{FE}, G } \atop {
\frac{Aq}{E},
\frac{Gq}{E},
\frac{A^2q^2}{CDEF} }} ; q \right) \cr
&\qquad  + {{(Aq,G,
\frac{Aq}{CE},
\frac{Aq}{DE},
\frac{Aq}{FE},
\frac{A^2q^2}{CDFG} )_\infty } \over {(
\frac{Aq}{C},
\frac{Aq}{D},
\frac{Aq}{E},
\frac{Aq}{F},
\frac{G}{E},
\frac{A^2q^2}{CDEFG} )_\infty }} {_4\phi_3} \left( {{
\frac{Aq}{CG},
\frac{Aq}{DG},
\frac{Aq}{FG}, E } \atop {
\frac{Aq}{G},
\frac{Eq}{G},
\frac{A^2q^2}{CDFG} }} ; q \right) . \cr
}
$$
The above is equivalent to the standard three-term 
transformation formula [7, ($\caprom 3$.~36), p.~246].

\bigskip
\item{3.} Replace $E$ and $F$ by $eq^{-m}$ and $fq^m$
respectively and take the limit as $m \rightarrow \infty$. 
The general formula reduces to a relation connecting
five $\ephis$'s. For particular values of the parameters 
this relation reduces to some known transformation formulas.
\par

A first particular case of the above is obtained by 
eliminating $f$ with the help of the balance condition 
and then writing $e=1$, $B=\qmn$. The result reduces to 
Sear's transformation of terminating balanced ${_4\phi_3}$ 
series [7, ($\caprom 3$.~15), p.~242].

A second particular case is obtained as follows. After 
eliminating $f$ as before, we take the special value 
$e = \frac{A^2q^2}{BDGH}$.
The resulting relation connects the three $\ephis$'s which are 
complement to each other viz.,
$$
\eqalign{
& \w \left(A;B,C,D,G,H; \frac{A^2q^2}{BCDGH} \right), \cr
& \w \left(\frac{B^2}{A};B,\frac{BC}{A},\frac{BD}{A},\frac{BG}{A},
\frac{BH}{A}; \frac{A^2q^2}{BCDGH} \right) \qquad {\rm and} \cr
& \w \left(\frac{G^2}{A};G,\frac{GB}{A},\frac{GC}{A},\frac{GD}{A},
\frac{GH}{A}; \frac{A^2q^2}{BCDGH} \right) . \cr
}
$$
An application of Bailey's infinite product identity 
[7, Ex. 5.21, p.~138) enables us to write the above transformation 
formula in the form given by Bailey [3] (also refer to 
[7, Ex. 2.15, p.~52].

\item{4.} Write $G=gq^m$, $H=hq^{-m}$ in (6.5) and take limit 
as $m \rightarrow \infty$. This results in the three-term
$\ephis$ transformation (5.8). 
\par

\bigskip
\noindent{\bf 7. Quadratic identities}
\medskip

Any two linearly independent solutions $X_n^{(i)}$ and $X_n^{(j)}$ 
of the three-term recurrence (3.2) satisfy the formula 
$$
\lim\limits_{n \rightarrow \infty} {{\CW (X_n^{(i)},X_n^{(j)})}
\over
{b_1b_2\dots b_n}} = b_0 \CW (X_{-1}^{(i)},X_{-1}^{(j)})
\leqno{(7.1)}
$$ 
where the Cassorati determinant
$$
\CW (X_n^{(i)},X_n^{(j)}):=X_n^{(i)}X_{n+1}^{(j)}-X_{n+1}^{(i)}X_n
^{(j)}.
$$
Taking second order asymptotics of (7.1) we may derive 
quadratic identities involving $\ephis$'s. We demonstrate this
below
by considering pairs of solutions from the solutions 
$X_n^{(1), \frac{s\qnm}{h}}$,
$X_n^{(2), \frac{\qnp}{h}}$ and
$X_n^{(3), \frac{bs}{ah}\qnm}$.

If $n$, $k$ are non-negative integers and $\alpha$, $\beta$ are 
independent of $n$, then
$$
\eqalign{
(\alpha \qmn )_k & = (1 - a_1 \frac{\qn}{\alpha} 
+ a_2 \frac{\qtn}{\alpha^2} + O(q^{3n}))(-\alpha)^k
q^{k(k-1)/2-nk}, \cr
(\beta \qn )_k & = 1 - b_1 \beta{\qn}
+ b_2 \beta^2{\qtn} + O(q^{3n}), \cr
}
\leqno{(7.2)}
$$
where 
$$
\eqalign{
& a_1 = \qmkp \frac{(1-\qk)}{(1-q)}, \qquad
a_2 = \qmkp \qmkpt \frac{(1-\qkm)(1-\qk)}{(1-q)(1-q^2)}, \cr
& b_1 = \frac{(1-\qk)}{(1-q)}, \qquad
b_2 = \frac{q(1-\qkm)(1-\qk)}{(1-q)(1-q^2)}. \cr
}
$$
Using (7.2) we can work out the second order  large $n$ asymptotics

of the three solutions mentioned above. Starting from (3.3) 
we obtain 
$$
\eqalign{
& X_n^{(1),\frac{s}{h}\qnm} \approx W_1 - \frac{\qn}{h} T_1W_{1+}
+ \frac{\qn}{h} S_1W_1 + O(\qtn) \cr
& W_1 = \w \left(a;b,c,d,e,f; \frac{s}{aq} \right) \cr
& W_{1+}=\w \left(aq^2;bq,cq,dq,eq,fq; \frac{s}{aq^2} \right) \cr
& S_1 = \frac{1}{(1-q)} \left( 
\frac{s}{q}
+ \frac{s}{aq}
+ \frac{aq}{b}
+ \frac{aq}{c}
+ \frac{aq}{d}
+ \frac{aq}{e}
+ \frac{aq}{f} - aq \right) \cr
& T_1 = {{\frac{s}{aq}
(1 - \frac{s}{aq^2})
(1-aq)
(1-aq^2)
(1-b)
(1-c)
(1-d)
(1-e)
(1-f) } \over {
(1-q)
(1 - \frac{aq}{b})
(1 - \frac{aq}{c})
(1 - \frac{aq}{d})
(1 - \frac{aq}{e})
(1 - \frac{aq}{f}) }}.  \cr
}
\leqno{(7.3)}
$$
In the above calculation we were surprised to find that not only
the first but also the second order asymptotics were given in 
terms of very well poised $\ephis$'s. We have not checked if this
continues to hold in higher orders. 

From (3.8) we obtain
$$
\eqalign{
& X_n^{(3), \frac{bs}{ah}\qnm} \approx L_3 \left(W_3 -
\frac{\qn}{h} 
T_3W_{3+} + \frac{\qn}{h}S_3W_3 \right)+ O (\qtn), \cr
& W_3 = \w \left( \frac{b^2}{a}; b, 
\frac{bc}{a},
\frac{bd}{a},
\frac{be}{a},
\frac{bf}{a};
\frac{s}{aq} \right) \cr
& W_{3+} = \w \left( \frac{b^2q^2}{a}; bq, 
\frac{bcq}{a},
\frac{bdq}{a},
\frac{beq}{a},
\frac{bfq}{a};
\frac{s}{aq^2} \right) \cr
& S_3 = (1-q)^{-1} \left(
\frac{aq}{c} +
\frac{aq}{d} +
\frac{aq}{e} +
\frac{aq}{f} +
\frac{bs}{aq} +
q +
\frac{s}{bq} - bq \right) \cr
& T_3 = {{\frac{s}{bq} 
(1 - \frac{s}{aq^2})
(1 - \frac{b^2q}{a})
(1 - \frac{b^2q^2}{a})
(1 - \frac{bc}{a})
(1 - \frac{bd}{a})
(1 - \frac{be}{a})
(1 - \frac{bf}{a})(1-b)} \over {(1-q)
(1 - \frac{bq}{c})
(1 - \frac{bq}{d})
(1 - \frac{bq}{e})
(1 - \frac{bq}{f})
(1 - \frac{bq}{a}) }} \cr
& L_3 = {{(
\frac{s}{aq},
\frac{bq}{c},
\frac{bq}{d},
\frac{bq}{e},
\frac{bq}{f},
\frac{bq}{a},
\frac{bhq}{s},
\frac{s}{bh} )_\infty } \over { (
\frac{b^2q}{a},
\frac{bc}{a},
\frac{bd}{a},
\frac{be}{a},
\frac{bf}{a}, b,
\frac{bh}{a},
\frac{aq}{bh} )_\infty }} . \cr
}
\leqno{(7.4)}
$$
From (3.6) we have
$$
\eqalign{
& X_n^{(2), \frac{\qnp}{h}} \approx W_2 - \frac{\qn}{h}T_2W_{2+}
+\frac{q^2}{h} S_2 W_2
+ O (\qtn) \cr
& W_2 = \w \left( 
\frac{q}{a};
\frac{q}{b},
\frac{q}{c},
\frac{q}{d},
\frac{q}{e},
\frac{q}{f};
\frac{aq^2}{s} \right) \cr
& W_{2+} = \w \left( 
\frac{q^3}{a};
\frac{q^2}{b},
\frac{q^2}{c},
\frac{q^2}{d},
\frac{q^2}{e},
\frac{q^2}{f};
\frac{aq}{s} \right) \cr
& S_2 = \frac{1}{(1-q)} \left(q+a+
\frac{bs}{aq} +
\frac{cs}{aq} +
\frac{ds}{aq} +
\frac{es}{aq} +
\frac{fs}{aq} -
\frac{s}{a} \right) \cr
& T_2 = {{a
(1 - \frac{aq}{s})
(1 - \frac{q^2}{a})
(1 - \frac{q^3}{a})
(1 - \frac{q}{b})
(1 - \frac{q}{c})
(1 - \frac{q}{d})
(1 - \frac{q}{e})
(1 - \frac{q}{f}) } \over {
(1 - q)
(1 - \frac{bq}{a})
(1 - \frac{cq}{a})
(1 - \frac{dq}{a})
(1 - \frac{eq}{a})
(1 - \frac{fq}{a}) }}.  \cr
}
\leqno{(7.5)}
$$
Looking at the convergence conditions of the different 
$\ephis$ series in (7.3), (7.4) and (7.5) we find that the 
arguments of series in (7.5) are not compatible with those of 
(7.3) and (7.4). The difficulty can be overcome by applying
the transformation [7, ($\caprom 3$.~39), p.~247]
to the $\Phi$ in $X_n^{(2), \frac{\qnp}{h}}$ before working out 
its asymptotics. We then obtain
$$
\eqalign{
& X_n^{(2), \frac{\qnp}{h}} \approx L_2'\left(W_2' - 
\frac{\qn}{h} T_2' W_{2+}' + \frac{\qn}{h} S_2' W_2' \right) + O 
(\qtn), \cr
& W_2' = \w \left(
\frac{cde}{a^2};
\frac{ce}{a},
\frac{cd}{a},
\frac{de}{a},
\frac{q}{f},
\frac{q}{b};
\frac{bf}{a} \right) \cr
& W_{2+}' = \w \left(
\frac{cdeq^2}{a^2};
\frac{ceq}{a},
\frac{cdq}{a},
\frac{deq}{a},
\frac{q^2}{f},
\frac{q^2}{b};
\frac{bf}{aq} \right) \cr
& S_2' = \frac{1}{(1-q)} \left(q+
\frac{cs}{aq} +
\frac{ds}{aq} +
\frac{es}{aq} +
\frac{a^2q}{cde} +
\frac{aq}{f} +
\frac{aq}{b} -
\frac{aq^2}{bf} \right) \cr
& T_2' = {{
\frac{bfs}{aq^2}
(1- \frac{bf}{aq})
(1- \frac{cdeq}{a^2})
(1- \frac{cdeq^2}{a^2})
(1- \frac{cd}{a})
(1- \frac{ce}{a})
(1- \frac{de}{a})
(1- \frac{q}{f})
(1- \frac{q}{b}) } \over {
(1-q)
(1 - \frac{cq}{a})
(1 - \frac{dq}{a})
(1 - \frac{eq}{a})
(1 - \frac{cdef}{a^2})
(1 - \frac{bcde}{a^2}) }}\cr
& L_2' = {{(
\frac{q^2}{a},
\frac{cdef}{a^2},
\frac{bcde}{a^2},
\frac{bf}{a}  )_\infty } \over {(
\frac{cdeq}{a^2},
\frac{fq}{a},
\frac{bq}{a},
\frac{aq^2}{s} )_\infty }} . \cr
}
\leqno{(7.6)}
$$
We now apply (7.1) to pairs of solutions from 
$X_n^{(1), \frac{s}{h}\qnm}$,
$X_n^{(2), \frac{\qnp}{h}}$ and
$X_n^{(3), \frac{bs}{ah}\qnm}$. In the particular case 
$h=1$, the right side of (7.1) can be evaluated in terms of
infinite 
products. Thus we obtain the three identities
$$
\leqalignno{
& L_2' (1-q) [(S_1 - S_2')W_1W_2' -T_1W_{1+} W_2' + T_2' W_1 
W_{2+}'] \cr
&\qquad = -q 
\left(1 - \frac{a}{b}\right) \left(1 - \frac{a}{c}\right)
\left(1 - \frac{a}{d}\right) \left(1 - \frac{a}{e}\right)
\left(1 - \frac{a}{f}\right)
\frac{(aq)_\infty}{(\frac{a}{q})_\infty} \cr 
&\qquad\qquad  \times \left[ 1 + {{(
\frac{a}{q}, \frac{q^2}{a}, b, \frac{q}{b}, c, \frac{q}{c},
d, \frac{q}{d}, e, \frac{q}{e},
f, \frac{q}{f},\frac{s}{q^2}, \frac{q^3}{s} )_\infty } \over {(a,
\frac{q}{a}, \frac{a}{b}, \frac{bq}{a}, \frac{a}{c},
\frac{cq}{a}, \frac{a}{d}, \frac{dq}{a}, \frac{a}{e},
\frac{eq}{a}, \frac{a}{f}, \frac{fq}{a}, \frac{s}{aq},
\frac{aq^2}{s} )_\infty }} \right] ,
& {(7.7)} \cr
& L_3'(1-q) [ (S_1 - S_3)W_1W_3 - T_1W_{1+}W_3 + T_3W_1W_{3+} ] \cr
&\qquad = - \frac{aq^2}{s} {{(aq, 
\frac{s}{bq}, \frac{bq^2}{s}, \frac{aq}{de}, \frac{aq}{ce},
\frac{aq}{cd}, \frac{aq}{cf}, \frac{aq}{df},
\frac{aq}{ef} )_\infty } \over {(
\frac{aq}{b}, \frac{aq}{c}, \frac{aq}{d}, \frac{aq}{e},
\frac{aq}{f}, \frac{bc}{a}, \frac{bd}{a}, \frac{be}{a},
\frac{bf}{a} )_\infty }} , 
& {(7.8)} \cr
& L_2' L_3' (1-q) [(S_2' - S_3)W_2'W_3 - T_2' W_{2+}'W_3 + 
T_3W_2'W_{3+}] \cr
&\qquad = - \frac{aq^2}{s} {{(
\frac{q^2}{a}, \frac{q}{c}, \frac{q}{d}, \frac{q}{e},
\frac{q}{f}, \frac{cdef}{a^2}, \frac{a^2q}{cdef}, \frac{s}{bq},
\frac{bq^2}{s} )_\infty } \over {(b,
\frac{b}{a}, \frac{bq}{a}, \frac{cq}{a}, \frac{dq}{a},
\frac{eq}{a}, \frac{fq}{a}, \frac{aq}{b},
\frac{aq^2}{s} )_\infty }} . 
& {(7.9)} \cr
}
$$
where $L_3'$ is the $L_3$ of (7.4) with $h=1$.

The identities (7.7) and (7.8) are easy to derive. However, in 
the derivation of (7.9) we had to make use of the $\Phi$ 
transformation formula (6.5) for special values 
$$
(A,B,C,D,E,F,G,H) \rightarrow \left(a,b,c,d,e,f,
\frac{s}{q},1\right)
$$
and also Slater's infinite product identity. 

It should be noted these identities are not independent. 
Starting from any two of the three identities (7.7), (7.8) and
(7.9)
we can deduce the third one. For example, if we multiply 
(7.7) by $\frac{W_3}{L_2'(1-q)}$
and (7.8) by $\frac{W_3}{L_3(1-q)}$ and subtract we obtain 
the value of
$$
(S_2' - S_3)W_2'W_3 - T_2' W_{2+}' W_3 + T_3 W_2' W_{3+}
$$
in terms of infinite products and the $W_1$, $W_2'$,
$W_3$. If we now apply the three-term $\ephis$ transformation 
formula  (5.8) and Slater's infinite product identity we can 
write down the result in the form of identity (7.9). 

We do not know where (7.7), (7.8) and (7.9) fit into the
general scheme of $q$-series identities. In order to
understand this better we now calculate a limiting case
of (7.8) where $f=\qmm$ and $m \rightarrow \infty$.

If $\left\vert \frac{aq}{de} \right\vert <1$ then, using
[7, ($\caprom 3$. 36), p.~246], we find that
$$
\lim_{m \rightarrow \infty} \w
\left(a;b,c,d,e,\qmm; \frac{a^2q^{m+2}}{bcde} \right)
= {{\left(aq, \frac{aq}{de} \right)_\infty}
\over {\left(\frac{aq}{d} , \frac{aq}{e} \right)_\infty
}} {}_3\phi_2 \left( {{\frac{aq}{bc} ,d,e } \atop
{\frac{aq}{b} , \frac{aq}{c} }} ; \frac{aq}{de} \right) ,
\leqno{(7.10)}
$$
$$
\lim_{m \rightarrow \infty} \w
\left(\frac{b^2}{a};b,\frac{bc}{a},\frac{bd}{a},
\frac{be}{a},\frac{b\qmm}{a}; \frac{a^3q^{m+3}}{bcde} \right)
= {{\left(\frac{b^2q}{a}, \frac{aq}{de} \right)_\infty}
\over {\left(\frac{bq}{d} , \frac{bq}{e} \right)_\infty
}} {}_3\phi_2 \left( {{\frac{q}{c} ,\frac{bd}{a},\frac{be}{a} 
} \atop {\frac{bq}{a} , \frac{bq}{c} }} ; \frac{aq}{de} 
\right) .
\leqno{(7.11)}
$$

This gives the $f = \qmm$, $m \rightarrow \infty$ limit of the
$W_1$ and $W_2$ in (7.8). We omit further details and
now state the final result.

 The $f = \qmm$, $m
\rightarrow \infty$ limit of (7.8) yields
$$
\eqalign{
& {}_3\phi_2 
\left( {{\frac{aq}{bc},d,e} \atop
{\frac{aq}{b}, \frac{aq}{c} }} ; \frac{aq}{de} \right)
{}_3\phi_2 
\left( {{\frac{q}{c}, \frac{bd}{a},
\frac{be}{a} } \atop {\frac{bq}{a}, \frac{bq}{c} }} ;
\frac{aq}{de} \right) \cr
& \qquad + 
\frac{aq}{cde} {{ (1-c)(1-d)(1-e)} \over {
\left( 1 - \frac{b}{a} \right)
\left( 1 - \frac{aq}{b} \right)
\left( 1 - \frac{aq}{c} \right) }} 
\rphit 
\left( {{\frac{aq}{bc},dq,eq} \atop
{\frac{aq^2}{b}, \frac{aq^2}{c} }} ; \frac{aq}{de} \right)
\rphit
\left( {{\frac{a}{c}, \frac{bd}{a},
\frac{be}{a} } \atop {\frac{bq}{a}, \frac{bq}{c} }} ;
\frac{aq}{de} \right) \cr
& \qquad - 
\frac{aq}{cde} {{ \left(1-\frac{bc}{a}\right)\left(1-
\frac{bd}{c}\right)\left(1-\frac{be}{a}\right)} \over {
\left( 1 - \frac{b}{a} \right)
\left( 1 - \frac{bq}{a} \right)
\left( 1 - \frac{bq}{c} \right) }} 
\rphit
\left( {{\frac{aq}{bc},d,e} \atop
{\frac{aq}{b}, \frac{aq}{c} }} ; \frac{aq}{de} \right)
{}_3\phi_2 
\left( {{\frac{q}{c}, \frac{bdq}{a},
\frac{beq}{a} } \atop {\frac{bq^2}{a}, \frac{bq^2}{c} }} ;
\frac{aq}{de} \right) \cr
& \qquad = {{\left( \frac{aq}{ce}, \frac{aq}{cd}, bq
\right)_\infty } \over { \left( \frac{aq}{c},
\frac{aq}{de} , \frac{bq}{c} \right)_\infty }} . \cr
}
\leqno{(7.12)}
$$
A further limit $c \rightarrow \infty$ is easily
calculated to give
$$
\eqalign{
& \tphio \left( {{d,e} \atop {\frac{aq}{b} }} ;
\frac{aq}{de} \right) \tphio \left( {{\frac{bd}{a} ,
\frac{be}{a} } \atop { \frac{bq}{a} }} ; \frac{aq}{de}
\right) \cr
& \qquad - 
\frac{aq}{de} {{ (1-c)(1-d)} \over {
\left( 1 - \frac{b}{a} \right)
\left( 1 - \frac{aq}{b} \right)
}} 
\tphio 
\left( {{dq,eq} \atop
{\frac{aq^2}{b} }} ; \frac{aq}{de} \right)
\tphio
\left( {{\frac{bd}{a},
\frac{be}{a} } \atop {\frac{bq}{a}, \frac{bq}{c} }} ;
\frac{aq}{de} \right) \cr
& \qquad + 
\frac{bq}{de} {{ \left(1-
\frac{bd}{a}\right)\left(1-\frac{be}{a}\right)} \over {
\left( 1 - \frac{b}{a} \right)
\left( 1 - \frac{bq}{a} \right)
}} 
\tphio
\left( {{d,e} \atop
{\frac{aq}{b} }} ; \frac{aq}{de} \right)
\tphio
\left( {{\frac{bdq}{a},
\frac{beq}{a} } \atop {\frac{bq^2}{a} }} ;
\frac{aq}{de} \right) \cr
& \qquad = {{(bq)_\infty} \over {\left( \frac{aq}{de}
\right)_\infty}} . \cr
}
\leqno{(7.13)}
$$

Even at the $\tphio$ level, this is not an identity we
are familiar with. However, if we now put $d=1$ it
reduces to
$$
\tphio \left( {{\frac{b}{a}, \frac{be}{a}} \atop {
\frac{bq}{a}}} ; \frac{aq}{e} \right) + \frac{bq}{de}{{\left(1-
\frac{be}{a} \right)}
\over {\left(1-\frac{bq}{a} \right)}} \tphio \left( {{\frac{bq}{a},
\frac{beq}{a}} \atop {\frac{bq^2}{a}}}; \frac{aq}{e}
\right) = {{(bq)_\infty} \over {\left(\frac{aq}{e}
\right)_\infty }}.
$$
This may now be recognized as the contiguous relation
[16]
$$
\tphio \left( {{A,B} \atop C};z\right) + {{Az(1-B)} \over
{(1-C)}} \tphio \left(  {{Aq, Bq} \atop {Cq}}; z \right)
= \tphio  \left( {{Aq, B} \atop {C}}; z \right)
$$
for the special case $C=Aq$ when
$$
\tphio \left( {{Aq, B} \atop {Aq}}; z \right) =
{}_1\phi_0 \left( {B \atop {-}} ; z \right) =
{{(Bz)_\infty} \over {(z)_\infty}} .
$$

\vglue 1 truein

\frenchspacing
\centerline{\bf References}

\item{1.\ }
R. Askey, Ramanujan and hypergeometric and basic
hypergeometric series, Ramanujan International Symposium on
Analysis, ed. N.~K.~Thakare, 1989, MacMillan India, Delhi,
1--83.
\medskip
\item{2.\ }
R. Askey and J. A. Wilson, Some basic hypergeometric
polynomials that generalize Jacobi polynomials, {\it Memoirs
Amer. Math. Soc.\/} {\bf 319} (1985), 1--55.
\medskip
\item{3.\ }
W. N. Bailey, Series of hypergeometric types which are
infinite in both directions, {\it Quart. J. Math.
(Oxford)\/} {\bf 7} (1936), 105--115.
\medskip
\item{4.\ }
W.~N.~Bailey, Generalized Hypergeometric series, Cambridge
University Press, Cambridge, reprinted by Hafner, New York,
1972. 
\medskip
\item{5.\ }
B.~C.~Berndt, R.~L.~Lamphere and B.~M.~Wilson, Chapter 12 of
Ramanujan's second notebook: Continued fractions, {\it Rocky
Mountain J. Math\/} {\bf 15} (1985), 235--310.
\medskip
\item{6.\ }
T.~S.~Chihara, An Introduction to Orthogonal Polynomials,
Gordon and Breach, New York, 1978.
\medskip
\item{7.\ }
G.~Gasper and M.~Rahman, Basic Hypergeometric Series,
Cambridge University Press, Cambridge, 1990.
\medskip
\item{8.\ }
W.~Gautschi, Computational aspects of three-term recurrence
relations, {\it SIAM Rev.\/} {\bf 9} (1967), 24--82.
\medskip
\item{9.\ }
D.~P.~Gupta and D.~R.~Masson,
Exceptional $q$-Askey-Wilson polynomials and continued
fractions, {\it Proc. Amer. Math. Soc.\/} {\bf 112} (1991),
717--727.
\medskip
\item{10.\ }
D.~P.~Gupta and D.~R.~Masson,
Watson's basic analogue of Ramanujan's Entry 40 and its
generalization, {\it SIAM J. Math. Anal.\/} {\bf 25} (1994),
429--440.
\medskip
\item{11.\ }
D.~P.~Gupta and D.~R.~Masson,
Solutions to the associated $q$-Askey-Wilson polynomial
recurrence relation, {\it Approximation and Computation\/},
Birkh\"auser, Boston, ed. R.~V.~M.~Zahar (1994), 273--284.
\medskip
\item{12.\ }
D.~P.~Gupta and D.~R.~Masson,
Contiguous relations, continued fractions and orthogonality:
an $\ephis$ model, {\it J. of Comp. and Appl. Math.\/} {\bf
65} (1995), to appear.
\medskip
\item{13.\ }
D.~P.~Gupta, M.~E.~H.~Ismail and D.~R.~Masson,
Associated continuous Hahn polynomials, {\it Canad. J.
Math.\/} {\bf 43} (1991), 1263--1280.
\medskip
\item{14.\ }
D.~P.~Gupta, M.~E.~H.~Ismail and D.~R.~Masson,
Contiguous relations, basic hypergeometric functions and
orthogonal polynomials $\caprom 2$, associated big
$q$-Jacobi polynomials, {\it J. Math. Anal. Appl.\/} {\bf
171} (1992), 477--497.
\medskip
\item{15.\ }
D.~P.~Gupta, M.~E.~H.~Ismail and D.~R.~Masson,
Contiguous relations, basic hypergeometric functions and
orthogonal polynomials $\caprom 3$, associated 
continuous dual $q$-Hahn polynomials, to appear in
{\it J. of Comp. and App. Math.\/}.
\medskip
\item{16.\ }
M.~E.~H.~Ismail and C.~Libis,
Contiguous relations, basic hypergeometric functions and
orthogonal polynomials,  
{\it J. Math. Anal. Appl.\/} {\bf 141} (1989), 349--372.
\medskip
\item{17.\ }
M.~E.~H.~Ismail and D.~R.~Masson,
Generalized orthogonality and continued fractions, {\it J.
Approx. Theory\/} {\bf 83} (1995), 1--40.
\medskip
\item{18.\ }
M.~E.~H.~Ismail and M.~Rahman,
Associated Askey-Wilson polynomials, {\it Trans. Amer. Math.
Soc.\/} {\bf 328} (1991), 201--239.
\medskip
\item{19.\ }
W.~B.~Jones and W.~J.~Thron, Continued Fractions: Analytic
Theory and Applications, Addison-Wesley, Reading, Mass.,
1980.
\medskip
\item{20.\ }
R. Koekoek and R.~F.~Swarttouw, The Askey-scheme of
hypergeometric orthogonal polynomials and its $q$-analogue,
Reports of the Faculty of Technical Mathematics and
Informatics no.~94-05, Delft, 1994.
\medskip
\item{21.\ }
D.~R.~Masson,
Some continued fractions of Ramanujan and Meixner-Pollaczek
polynomials, {\it Canad. Math. Bull.\/} {\bf 32} (1989),
177--181.
\medskip
\item{22.\ }
D.~R.~Masson,
Wilson polynomials and some continued fractions of
Ramanujan, {\it Rocky Mountain J. Math.} {\bf 21} (1991),
489--499.
\medskip
\item{23.\ } 
D.~R.~Masson, 
Associated Wilson polynomials, {\it Constr. Approx.\/} {\bf
7} (1991), 521--534.
\medskip
\item{24.\ }
D.~R.~Masson,
A generalization of Ramanujan's best theorem on continued
fractions, {\it Canad. Math. Reports of the Acad. of Sci.\/}
{\bf 13} (1991) 167--172.
\medskip
\item{25.\ }
D.~R.~Masson,
The last of the hypergeometric continued fractions, in 
Proceedings of the International Conference on Mathematical
Analysis and Signal Processing, Cairo, Egypt, January
3--9, 1994, to appear.
\medskip
\item{26.\ }
M.~Rahman and S.~K.~Suslov,
Classical biorthogonal rational functions, in Methods
of Approximation Theory in Complex Analysis and Mathematical
Physics $\caprom 4$, A.~A.~Gonchar and E.~B.~Saff, Eds.,
Lecture Notes in Mathematics 1550, Springer-Verlag, Berlin,
1993, 131--146.
\medskip
\item{27.\ }
L.~J.~Slater,
A note on equivalent product theorems, {\it Math. Gazette\/}
{\bf 38} (1954) 127--128.
\medskip
\item{28.\ }
G.~N.~Watson, 
Ramanujan's continued fraction, {\it Proc. Cambridge Philos.
Soc.\/} {\bf 31} (1935), 7--17.
\medskip
\item{29.\ }
J.~A.~Wilson, 
Hypergeometric Series, Recurrence relations, and Some New
Orthogonal Functions, Ph.D. Thesis, University of Wisconsin,
Madison, WI, 1978.
\medskip
\item{30.\ }
J.~A.~Wilson, 
Orthogonal functions from Gramm determinants, {\it SIAM J.
Math. Anal.\/} {\bf 22} (1991), 1147--1155.

\bye